\documentclass[reqno]{amsart}

\usepackage[T1]{fontenc}
\usepackage{lmodern}
\usepackage{amsmath,amssymb,amsthm,mathtools}
\usepackage{graphicx}
\usepackage{enumitem}
\usepackage{microtype}
\usepackage{cite}
\usepackage{hyperref}
\usepackage[noabbrev]{cleveref}

\numberwithin{equation}{section}
\allowdisplaybreaks

\theoremstyle{plain}
\newtheorem{theorem}{Theorem}[section]
\newtheorem{lemma}[theorem]{Lemma}
\newtheorem{proposition}[theorem]{Proposition}
\newtheorem{corollary}[theorem]{Corollary}
\theoremstyle{definition}
\newtheorem{definition}[theorem]{Definition}

\theoremstyle{remark}
\newtheorem{remark}[theorem]{Remark}
\newtheorem*{remark*}{Remark}

\crefname{theorem}{Theorem}{Theorems}
\Crefname{theorem}{Theorem}{Theorems}
\crefname{lemma}{Lemma}{Lemmas}
\Crefname{lemma}{Lemma}{Lemmas}
\crefname{proposition}{Proposition}{Propositions}
\Crefname{proposition}{Proposition}{Propositions}
\crefname{corollary}{Corollary}{Corollaries}
\Crefname{corollary}{Corollary}{Corollaries}
\crefname{definition}{Definition}{Definitions}
\Crefname{definition}{Definition}{Definitions}
\crefname{remark}{Remark}{Remarks}
\Crefname{remark}{Remark}{Remarks}
\crefname{section}{Section}{Sections}
\Crefname{section}{Section}{Sections}

\crefformat{section}{Section~#2#1#3}
\Crefformat{section}{Section~#2#1#3}
\crefformat{subsection}{Section~#2#1#3}
\Crefformat{subsection}{Section~#2#1#3}
\crefformat{theorem}{Theorem~#2#1#3}
\Crefformat{theorem}{Theorem~#2#1#3}
\crefformat{lemma}{Lemma~#2#1#3}
\Crefformat{lemma}{Lemma~#2#1#3}
\crefformat{proposition}{Proposition~#2#1#3}
\Crefformat{proposition}{Proposition~#2#1#3}
\crefformat{corollary}{Corollary~#2#1#3}
\Crefformat{corollary}{Corollary~#2#1#3}
\crefformat{definition}{Definition~#2#1#3}
\Crefformat{definition}{Definition~#2#1#3}
\crefformat{remark}{Remark~#2#1#3}
\Crefformat{remark}{Remark~#2#1#3}

\newcommand{\R}{\mathbb R}

\DeclareMathOperator{\supp}{supp}

\newcommand{\dd}{\,d}

\newcommand{\barredintegral}[2]{%
  \begingroup
  \setbox0=\hbox{$\mathsurround=0pt#1\int$}%
  \rlap{\hbox to\wd0{\hss\scalebox{1.3}[1]{$\mathsurround=0pt#1-$}\hss}}%
  \endgroup}

\crefname{assumption}{Assumption}{Assumptions}
\Crefname{assumption}{Assumption}{Assumptions}
\crefname{example}{Example}{Examples}
\Crefname{example}{Example}{Examples}
\crefformat{assumption}{Assumption~#2#1#3}
\Crefformat{assumption}{Assumption~#2#1#3}
\crefformat{example}{Example~#2#1#3}
\Crefformat{example}{Example~#2#1#3}

\AddToHook{cmd/maketitle/after}{\enlargethispage{3pt}}

\makeatletter
\@namedef{subjclassname@2020}{%
	\textup{2020} Mathematics Subject Classification}
\makeatother

\newcommand{\Lp}{\mathcal L_{s,p}}
\newcommand{\E}{\mathcal E}
\newcommand{\pv}{\operatorname{P.V.}}
\newcommand{\loc}{\mathrm{loc}}

\newcommand{\pair}[2]{\left\langle#1,#2\right\rangle}

\hypersetup{
  pdftitle={Weak--viscosity equivalence for the fractional Trudinger equation},
  pdfauthor={Qi Xue and Chao Zhang},
  pdfsubject={Weak--viscosity equivalence for the fractional Trudinger equation},
  pdfkeywords={Fractional Trudinger equation; weak--viscosity equivalence; nonlocal parabolic equations; multiplicative regularization.}
}

\begin{document}

\title[Weak--viscosity equivalence for the fractional Trudinger equation]
{Weak--viscosity equivalence for the fractional Trudinger equation}

\author{Qi Xue}
\address{Qi Xue \hfill\break School of Mathematics, Harbin Institute of Technology,
	Harbin 150001, P.R. China}
\email{18332276756@163.com}

\author{Chao Zhang}
\address{Chao Zhang \hfill\break School of Mathematics and Institute for Advanced Study in
	Mathematics, Harbin Institute of Technology, Harbin 150001, P.R. China}
\email{czhangmath@hit.edu.cn}
\thanks{Corresponding author: Chao Zhang.}
\date{}


\begin{abstract}
We prove weak--viscosity equivalence for the fractional Trudinger equation
\(\partial_t(|u|^{p-2}u)+\mathcal L_{s,p}u=0\),
where \(1<p<\infty\) and \(0<s<1\).
The functions are assumed to be continuous in the equation domain
and continuous in time with values in the weighted tail space
\(L^{p-1}_{sp}(\R^N)\). We treat both functions that are strictly
positive in the equation domain, with signed exterior values,
and globally nonnegative functions with interior zeros.
Exterior values need not be bounded or spatially continuous.
For \(1<p<2\), the viscosity formulation tests the evolution
variable \(|u|^{p-2}u\), retaining the nonlocal interaction at zeros.
The viscosity-to-weak implication follows from localized
multiplicative regularization and energy estimates independent
of a positive lower bound. The reverse implication uses comparison
with smooth functions, justified by a relative-energy identity.
A zero-slice property allows the weak identities on positive time
intervals to be joined across the zero set.
\end{abstract}

\maketitle

\section{Introduction and main results}
\label{sec:intro}

Let \(\Omega\subset\R^N\) be open, \(I=(0,T)\), and
\(Q=\Omega\times I\). Set
\begin{equation}
 J_p(a)=|a|^{p-2}a,\qquad
 \Lp u(x,t)=2\pv\int_{\R^N}
 \frac{J_p(u(x,t)-u(y,t))}{|x-y|^{N+sp}}\dd y.
 \label{eq:operator}
\end{equation}
We consider
\begin{equation}
 \partial_tJ_p(u)+\Lp u=0\quad\hbox{in }Q.
 \label{eq:main}
\end{equation}
Throughout the paper, \(N\ge1\), \(1<p<\infty\), and \(0<s<1\),
unless a smaller parameter range is specified.
The normalization in \eqref{eq:operator} is consistent with the
energy form \eqref{eq:energy}.

\subsection{Main results}

Write
\[
 \varpi(x)=\frac{1}{1+|x|^{N+sp}},\qquad
 X:=L^{p-1}_{sp}(\R^N)=L^{p-1}(\R^N,\varpi(x)\dd x).
\]
All functions are measurable on \(\R^N\times I\). We assume
\begin{equation}
 u\in C(Q)\cap C_{\loc}(I;X).
 \label{eq:admissible-class}
\end{equation}
For every \(p>1\), the second condition means
\[
 \int_{\R^N}|u(x,t)-u(x,t_0)|^{p-1}\varpi(x)\dd x
 \longrightarrow0\qquad(t\to t_0).
\]
We write \(\|f\|_X=(\int |f|^{p-1}\varpi)^{1/(p-1)}\).
This is a norm for \(p\ge2\) and a quasi-norm for \(1<p<2\).
In the latter case the topology is induced by the metric
\(d_X(f,g)=\int|f-g|^{p-1}\varpi\).
The notation \(L^\infty(J;X)\) means that the tail moments
\(\int |u(x,t)|^{p-1}\varpi(x)\dd x\) are essentially bounded on \(J\).
We use the continuous \(X\)-valued representative at every time; outside \(\Omega\),
values on spatial null sets are immaterial. The first condition
implies local boundedness in \(Q\). No boundedness on the full
strip, and no pointwise continuity of the exterior values, is assumed.

For the positive case, we use the class
\begin{equation}
 u\in C(Q)\cap C_\loc(I;X),
 \qquad u>0\quad\hbox{in }Q.
 \label{eq:class}
\end{equation}
In particular, \(\inf_{Q'}u>0\) for every \(Q'\Subset Q\).
Neither a positive lower bound up to \(\partial Q\) nor a sign
condition on \((\R^N\setminus\Omega)\times I\) is assumed.
The notions of weak and viscosity solution are defined in
Section~\ref{sec:framework}.

\begin{theorem}
\label{thm:main}
Let \(1<p<\infty\), \(0<s<1\), and let \(u\) satisfy \eqref{eq:class}.
Then \(u\) is a local weak solution of \eqref{eq:main} in \(Q\) if and
only if it is a viscosity solution there.
Moreover, every viscosity supersolution, respectively subsolution,
satisfying \eqref{eq:class} belongs to
\(L^p_\loc(I;W^{s,p}_\loc(\Omega))\) and is a local weak
supersolution, respectively subsolution.
\end{theorem}

The viscosity formulation is given in Definition~\ref{def:viscosity}
for \(p\ge2\), and in Definition~\ref{def:viscosity-singular}
for \(1<p<2\). In a positive phase the latter is equivalent to
testing \(u\) itself, with the spatial admissibility condition
\eqref{eq:critical-admissibility} when necessary.
The one-sided assertion concerns the viscosity-to-weak direction;
the reverse implication uses the weak equation with equality.

Nonnegativity on the full spatial strip allows interior zeros.
The exterior values may still be unbounded and time-dependent.

\begin{theorem}
\label{thm:nonnegative}
Let \(1<p<\infty\), \(0<s<1\), let \(u\) satisfy
\eqref{eq:admissible-class}, and suppose that
\[
 u(\cdot,t)\ge0\quad\hbox{almost everywhere on }\R^N
 \quad\hbox{for every }t\in I.
\]
Then \(u\) is a local weak solution of \eqref{eq:main} in \(Q\)
if and only if it is a viscosity solution there.
Moreover, every viscosity subsolution in this class belongs to
\(L^p_\loc(I;W^{s,p}_\loc(\Omega))\) and is a local weak
subsolution.
\end{theorem}

For a sign-changing viscosity solution \(u\) in \eqref{eq:admissible-class},
Proposition~\ref{prop:nonnegative-subsolution} applied to \(u\) and \(-u\)
gives \(u\in L^p_\loc(I;W^{s,p}_\loc(\Omega))\). For \(p\ne2\), the
argument gives the weak equation for such \(u\) only on \(\{u\ne0\}\).

\subsection{Related results}

For the local Trudinger equation, Lindgren and Lindqvist established a
comparison result and a weak-to-viscosity implication
\cite{LL22}. Lindqvist, Parviainen, and Siltakoski proved equivalence for
strictly positive local solutions
\cite[Theorems~2.8 and~2.9]{LPS26}. For \(p\ge2\), their subsequent
preprint establishes equivalence for nonnegative solutions and proves
the weak-to-viscosity implication also for signed solutions
\cite[Theorems~1 and~2]{LPSnonnegative}.
Anttila, Lindqvist, and Parviainen established comparison under a
strict gap between the boundary values
and, for \(p\ge2\), proved uniqueness
for signed weak solutions with time-dependent \(C^2\) boundary data
\cite[Theorems~3.1 and~4.1]{ALP26}.

In the nonlocal elliptic setting, Korvenp\"a\"a, Kuusi, and Lindgren
proved equivalence of solution notions for fractional \(p\)-Laplace
type equations and developed the critical-point test class needed
when \(p\le2/(2-s)\) \cite{KKL19}. Related nonhomogeneous elliptic problems
were studied by Barrios and Medina \cite{BM}. Ghosh, Lakshmi, and Zhang
proved weak--viscosity equivalence for nonlocal elliptic double phase
problems with nonhomogeneous data \cite{GLZ26}. In the evolutionary
setting, Jesus, Sobral, and Urbano study the fractional
\(p\)-caloric equation \(u_t+\Lp u=0\) in the degenerate range.
Their weak-to-viscosity result assumes interior continuity and
time continuity with values in \(X\); the reverse implication requires
global boundedness and semicontinuity on the full spatial strip
\cite[Definition~2.2 and Theorems~3.3--3.4]{JSU26}.
For \eqref{eq:main}, the time derivative of \(J_p(u)\) requires
a separate treatment at zeros of \(u\) when \(1<p<2\).

For the evolution equation \(J_p(u_t)+\Lp u=0\), Hynd and
Lindgren related weak and viscosity solutions of an
initial-boundary value problem \cite[Theorem 1.3]{HL16}.
They also studied weak solutions and their large-time behavior for
the fractional Trudinger equation with zero exterior values
\cite[Section 6]{HL21}. For weak solutions of the latter equation,
Adimurthi proved local H\"older regularity for bounded signed
solutions \cite{Adimurthi25}. Prasad proved a weak Harnack
inequality for globally nonnegative weak supersolutions
\cite{Prasad24}, while Ciani and Nakamura obtained refined
Harnack estimates \cite{CN26}.

The change of variable \(w=J_p(u)\) writes \eqref{eq:main} as
\(w_t+\Lp(J_p^{-1}(w))=0\). Collier and Hauer study this type of
nonlinear diffusion in a semigroup framework. Their comparison
theorem allows signed, time-dependent exterior data for strong
distributional solutions, with both evolution variables in
\(W^{1,1}(0,T;L^1(\Omega))\)
\cite[Theorem~8.1]{CH22}.
Strunk established a comparison principle for nonnegative weak
solutions of fractional doubly nonlinear equations under the
assumption that at least one exterior datum is independent of
time \cite[Theorem 1.4]{Strunk26}.

\subsection{Outline of the proof}

Equation~\eqref{eq:main} is invariant under
positive scalar multiplication, but generally not under
the addition of constants when \(p\ne2\).
We therefore use a multiplicative regularization
confined to an interior region, leaving the exterior values
unchanged.
This introduces an error in the viscosity inequality, controlled
by translation continuity in \(X\). Thus the weak formulation
can be recovered without imposing spatial continuity or boundedness
on the exterior data.

In the range \(1<p\le2/(2-s)\), the spatial penalty must have
order greater than \(sp/(p-1)\). We keep the time penalty quadratic
and treat the two variables separately. At each point of the lower
envelope, a minimizing shift gives an upper support and hence a lower
bound for the truncated operators, uniform in the truncation radius.
At points of differentiability with zero spatial gradient, a flat lower
test for the attaining translate gives the differential inequality for
the envelope. The upper envelope satisfies the reverse inequalities.
The resulting energy estimates are independent of the positive
lower bound. For the converse, a relative-energy identity justifies
comparison between a weak solution and a smooth subsolution or
supersolution. Freezing the spatial profile of a viscosity test
preserves its critical-point admissibility on a backward cylinder.
Only this restricted comparison is needed; comparison between
arbitrary weak solutions is not asserted.

To treat nonnegative functions, we first obtain local energy bounds
by positive truncation, as in \cite[Section~4]{LPSnonnegative}.
An interior zero of a globally nonnegative viscosity supersolution
forces the entire spatial slice to vanish almost everywhere.
The weak identities on the positive time intervals can then be
summed across the zero set. In the reverse direction, the weak
Harnack inequality of \cite{Prasad24} and the weak equation give
the corresponding zero-slice property. For \(1<p<2\), the viscosity formulation
uses \(v=J_p(u)\). At a zero contact point, the nonlinear
numerator of the spatial integrand reduces to the negative
of the patched test for \(v\). This retains the exterior interaction despite the
singularity of \(J_p'\) at zero.

The paper is organized as follows. Section~\ref{sec:framework} gives the definitions and preliminary
estimates. Sections~\ref{sec:envelopes} and~\ref{sec:regularized}
prove the viscosity-to-weak implication in the positive case;
Section~\ref{sec:comparison} proves the converse.
Section~\ref{sec:nonnegative} treats nonnegative solutions, and
Section~\ref{sec:scope} discusses the scope of the assumptions.

\section{Preliminaries}
\label{sec:framework}

\subsection{Weak formulation and energy tests}

Set \(p'=p/(p-1)\) and \(\beta=J_p\). For functions for which
the integral is finite, define
\begin{equation}
 \E(w,\eta)=\iint_{\R^{2N}}
 \frac{J_p(w(x)-w(y))(\eta(x)-\eta(y))}{|x-y|^{N+sp}}
 \dd x\dd y.
 \label{eq:energy}
\end{equation}
For a bounded open set \(D\), put
\[
 [w]_{s,p;D}^p=\iint_{D\times D}
 \frac{|w(x)-w(y)|^p}{|x-y|^{N+sp}}\dd x\dd y,
 \qquad
 V(D)=\overline{C_c^\infty(D)}^{\,W^{s,p}(\R^N)}.
\]
Functions in \(V(D)\) vanish almost everywhere outside \(D\).
We write \([w]_{s,p}\) for the seminorm on \(\R^N\).

\begin{definition}[Weak solutions]
\label{def:weak}
A measurable function \(u:\R^N\times I\to\R\) with
\(u\in L^\infty_\loc(Q)\cap L^\infty_\loc(I;X)\) is a local weak solution of
\eqref{eq:main} in \(Q\) if
\(u\in L^p_\loc(I;W^{s,p}_\loc(\Omega))\) and
\begin{equation}
 -\iint_Q\beta(u)\eta_t\dd x\dd t
 +\int_I\E(u(\cdot,t),\eta(\cdot,t))\dd t=0
 \label{eq:weak}
\end{equation}
for every \(\eta\in C_c^\infty(Q)\). In the energy form, \(\eta\) is
extended by zero outside \(\Omega\). Under the same assumptions,
\(u\) is a local weak supersolution or subsolution if the
left-hand side of \eqref{eq:weak} is
nonnegative or nonpositive, respectively, for every nonnegative
\(\eta\in C_c^\infty(Q)\).
\end{definition}

The testing form is finite by H\"older's inequality near the test
support and by the tail bound for the separated interactions.
Although the equation is imposed only in \(Q\), the function is
specified on \(\R^N\times I\). No initial or boundary condition
is included in Definition~\ref{def:weak}.

\begin{lemma}
\label{lem:energy-testing}
Let \(D\Subset D'\Subset\Omega\) be bounded open sets and
\(J\Subset I\). Suppose that \(|U|\le M\) on \(D'\times J\),
\(U\in L^\infty(J;X)\), and \(U\in L^p(J;W^{s,p}(D'))\). The operator
\[
 \pair{\mathcal A_U(t)}{\eta}=\E(U(\cdot,t),\eta),
 \qquad \eta\in V(D),
\]
belongs to \(L^{p'}(J;V(D)^*)\), and for almost every \(t\),
\begin{equation}
 \|\mathcal A_U(t)\|_{V(D)^*}
 \le C\bigl([U(\cdot,t)]_{s,p;D'}^{p-1}+M^{p-1}
       +\|U(\cdot,t)\|_X^{p-1}\bigr).
 \label{eq:dual-operator-bound}
\end{equation}
Here \(C\) depends on \(N,s,p,D,D'\).
If \(U\) is a weak solution near \(\overline D\times\overline J\),
then \(\partial_t\beta(U)=-\mathcal A_U\) in this dual space.

More generally, if a weak subsolution or supersolution also has
\(\partial_t\beta(U)=a\in L^{p'}(J;V(D)^*)\), its inequality
holds in the form
\[
 \int_J\bigl(\pair{a}{\eta}+\E(U,\eta)\bigr)\dd t
 \ \le0\quad\hbox{or}\quad\ge0
\]
for nonnegative \(\eta\in L^p(J;V(D))\) supported in a fixed
compact spatial subset of \(D\).
\end{lemma}

\begin{proof}
On \(D'\times D'\), H\"older's inequality bounds the testing
form by
\([U]_{s,p;D'}^{p-1}[\eta]_{s,p}\).
In the remaining nonzero pairs, one point belongs to \(D\) and
the other lies outside \(D'\). Their separation is positive, so
the kernel is bounded by \(C\varpi(y)\) when \(x\in D\) and
\(y\notin D'\). This contribution is therefore bounded by
\(C(M^{p-1}+\|U(\cdot,t)\|_X^{p-1})\|\eta\|_{L^p(D)}\).
This proves \eqref{eq:dual-operator-bound}; the identity for weak
solutions follows from the distributional equation.

If the stated dual time derivative exists, integration by parts
first rewrites the weak inequality for smooth tests as a pairing
with \(a+\mathcal A_U\). A nonnegative energy test with compact
spatial support is approximated by nonnegative space-time
mollifications, after a temporal cutoff. These converge in
\(L^p(J;W^{s,p}(\R^N))\), and their supports remain inside \(D\)
for small mollification radii. Both dual pairings therefore
converge. Removing the temporal cutoff by the same norm convergence
proves the assertion.
\end{proof}

\begin{remark}
\label{rem:time-traces}
Definition~\ref{def:weak} uses compactly supported space-time
tests. In the equivalence theorems, \(u\in C(Q)\) is assumed.
Consequently \(u\) and \(\beta(u)\) have continuous local
\(L^r\) traces for every finite \(r\), on compact interior
cylinders. A temporal cutoff in \eqref{eq:weak} therefore gives
\[
 \left[\int_\Omega\beta(u)\eta\dd x\right]_{t_1}^{t_2}
 -\int_{t_1}^{t_2}\int_\Omega\beta(u)\eta_t\dd x\dd t
 +\int_{t_1}^{t_2}\E(u,\eta)\dd t=0
\]
for every \([t_1,t_2]\Subset I\) and every smooth test \(\eta\)
whose spatial support lies in a fixed compact subset of \(\Omega\).
The testing form is integrable by Lemma~\ref{lem:energy-testing}.
For \(D\Subset\Omega\), the identity extends
by density to tests
\(\eta\in L^p(t_1,t_2;V(D))
\cap W^{1,p}(t_1,t_2;L^p(D))\)
whose spatial supports lie in a fixed compact subset of \(D\).
The spatial term passes to the limit by
Lemma~\ref{lem:energy-testing}, while
\(\beta(u)\in C([t_1,t_2];L^{p'}(D))\)
justifies convergence of the time and endpoint terms.
This argument applies to every \(p>1\); it does not
require \(\beta\) to be Lipschitz at zero or the embedding
\(V(D)\hookrightarrow L^2(D)\).
Time continuity need not hold for one-sided weak solutions:
the spatially constant function \(u(x,t)=\mathbf1_{\{t>t_0\}}\)
is a weak supersolution with a jump at \(t_0\in I\).
This explains the additional time-derivative hypothesis for the
one-sided assertion of Lemma~\ref{lem:energy-testing}.
\end{remark}

\subsection{Viscosity tests and patched operators}

We use \(C^{2,1}\) for functions whose spatial derivatives up to
order two and first time derivative are continuous. In a viscosity
test, the solution is replaced only inside the contact ball;
its exterior values are retained.

Put \(p_*=2/(2-s)\). A \(C^2\) function \(f\), defined near
\(x_0\), is spatially admissible at \(x_0\) if \(p>p_*\), or
\(Df(x_0)\ne0\), or the following condition holds:
\begin{equation}
 \begin{gathered}
 Df(x_0)=0,\qquad \kappa>\frac{sp}{p-1},\\
 c|x-x_0|^{\kappa-1}\le |Df(x)|,\qquad
 |D^2f(x)|\le C|x-x_0|^{\kappa-2}
 \quad(0<|x-x_0|<r)
 \end{gathered}
 \label{eq:critical-admissibility}
\end{equation}
for some \(r,c,C>0\).
Only when \(p\le p_*\) and \(Df(x_0)=0\) is this extra
condition imposed. It is the local form of the isolated-critical-point
\(C^2_\kappa\) condition in \cite{KKL19}. In particular,
\(x_0\) is the only critical point in \(B_r(x_0)\).

\begin{definition}[Viscosity solutions]
\label{def:viscosity}
Let \(p\ge2\), \(0<s<1\), and let \(u\) satisfy
\eqref{eq:admissible-class}.
We call \(u\) a viscosity supersolution of \eqref{eq:main} in \(Q\)
if the following condition holds for every \(z_0=(x_0,t_0)\in Q\)
and every \(\rho,a>0\) such that
\[
 C^-_{\rho,a}(z_0):=
 \overline{B_\rho(x_0)}\times[t_0-a,t_0]\Subset Q.
\]
Let \(\phi\in C^{2,1}\) on a neighborhood of
\(C^-_{\rho,a}(z_0)\) satisfy
\[
 \phi(z_0)=u(z_0),\qquad
 \phi\le u\quad\hbox{on }C^-_{\rho,a}(z_0).
\]
For \((y,t)\in\R^N\times[t_0-a,t_0]\), define
\[
 \phi_{\rho,u}(y,t)=
 \begin{cases}
  \phi(y,t),&y\in B_\rho(x_0),\\
  u(y,t),&y\notin B_\rho(x_0).
 \end{cases}
\]
Then, for every such \(\phi\), we require
\begin{equation}
 \beta'(u(z_0))\phi_t(z_0)+\Lp\phi_{\rho,u}(z_0)\ge0.
 \label{eq:viscosity}
\end{equation}
A viscosity subsolution is defined by requiring
\(\phi(z_0)=u(z_0)\), \(\phi\ge u\) on \(C^-_{\rho,a}(z_0)\),
and the reverse inequality in \eqref{eq:viscosity}.
A viscosity solution is both a viscosity subsolution and a viscosity
supersolution. Here \(\beta'(r)=(p-1)|r|^{p-2}\) for \(p>2\),
and \(\beta'(r)=1\) for \(p=2\).
\end{definition}

The contact is only backward in time and need not be strict.
This backward-cylinder and exterior-patching convention is also used
for the fractional \(p\)-caloric equation in
\cite[Definition~2.2]{JSU26}. At the contact point, the spatial term is
\[
\begin{aligned}
 \Lp\phi_{\rho,u}(z_0)
 ={}&2\pv\int_{B_\rho(x_0)}
 \frac{J_p(\phi(z_0)-\phi(y,t_0))}{|x_0-y|^{N+sp}}\dd y\\
 &+2\int_{\R^N\setminus B_\rho(x_0)}
 \frac{J_p(u(z_0)-u(y,t_0))}{|x_0-y|^{N+sp}}\dd y.
\end{aligned}
\]
Thus the exterior values are evaluated at the same time \(t_0\);
no extension of \(\phi\) to the whole space is required.
The definition allows signed functions. Strict positivity is an
assumption of Theorem~\ref{thm:main}, not of the viscosity definition.
In a positive phase the time coefficient is \((p-1)u(z_0)^{p-2}\).

For \(1<p<2\), that coefficient is singular at zero. We instead
test the evolution variable \(v=\beta(u)\). Write
\(\gamma(r)=\operatorname{sgn}(r)|r|^{1/(p-1)}\), with
\(\gamma(0)=0\).

\begin{definition}
\label{def:viscosity-singular}
Let \(1<p<2\), and let \(u\) satisfy
\eqref{eq:admissible-class}. Set \(v=\beta(u)\).
We call \(u\) a viscosity supersolution if, whenever
\(\psi\in C^{2,1}\) touches \(v\) from below on a backward
cylinder \(C^-_{\rho,a}(z_0)\Subset Q\), the inequality
\begin{equation}
 \psi_t(z_0)+\Lp\bigl(\gamma(\psi_{\rho,v})\bigr)(z_0)\ge0
 \label{eq:singular-viscosity}
\end{equation}
holds for the following admissible tests.
If \(v(z_0)\ne0\), the function
\(\gamma(\psi(\cdot,t_0))\), which is smooth near \(x_0\),
must be spatially admissible at \(x_0\).
If \(v(z_0)=0\), every such \(C^{2,1}\) test is admitted,
including the constant zero test.
Here \(\psi_{\rho,v}=\psi\) in the ball and \(v\) outside it,
with the same patching convention as in
Definition~\ref{def:viscosity}.
Subsolutions are defined by upper contact and the reversed inequality;
solutions satisfy both inequalities.
\end{definition}

The exterior values of \(\gamma(\psi_{\rho,v})\) are \(u\).
At a zero contact value, \(J_p\circ\gamma=\operatorname{id}\) gives
\begin{equation}
 \Lp\bigl(\gamma(\psi_{\rho,v})\bigr)(z_0)
 =-2\pv\int_{\R^N}
       \frac{\psi_{\rho,v}(y,t_0)}{|x_0-y|^{N+sp}}\dd y.
 \label{eq:zero-transformed-operator}
\end{equation}
The local linear Taylor term cancels, and the remainder is
integrable because \(sp<2\). The exterior integral is finite
since \(|v|=|u|^{p-1}\). Thus the principal value in Definition~\ref{def:viscosity-singular}
is finite at zero contacts without a spatial flatness restriction.
The example in \eqref{eq:exterior-bump-residual} shows how the exterior
values enter the equation at such a contact.

\begin{lemma}
\label{lem:positive-tests}
Let \(1<p<2\) and \(u>0\) in \(Q\).
Definition~\ref{def:viscosity-singular} is equivalent to the
following formulation: test \(u\) by \(C^{2,1}\) functions
\(\phi\), require \(\phi(\cdot,t_0)\) to be spatially
admissible at each contact, and impose
\[
 (p-1)u(z_0)^{p-2}\phi_t(z_0)
       +\Lp\phi_{\rho,u}(z_0)\gtreqless0
\]
for lower and upper tests, respectively.
\end{lemma}

\begin{proof}
Near a positive contact, the substitutions \(\psi=\beta(\phi)\)
and \(\phi=\gamma(\psi)\) preserve the order of contact and
satisfy \(\psi_t=\beta'(u)\phi_t\).
If needed, first shrink the contact ball and time interval so
that the tests remain positive. The compositions are then
\(C^{2,1}\), and the admission rule is precisely the one stated
for \(\phi\).
Restoring a larger original patch preserves the inequality:
for a lower test its values on the added annulus are at most
\(u\), which increases the operator; for an upper test the
operator decreases. This proves both implications for every
original patching radius.
\end{proof}

For \(p>p_*\), the following estimate ensures that the principal
value exists at every contact point, including spatial critical
points. Below this threshold, we use the flatness condition above.

\begin{lemma}
\label{lem:operator}
For \(p\ge2\),
\begin{equation}
 |J_p(a)-J_p(b)|\le C_p(|a|+|b|)^{p-2}|a-b|.
 \label{eq:jp-estimate}
\end{equation}
For \(1<p<2\),
\begin{equation}
 |J_p(a)-J_p(b)|\le C_p|a-b|^{p-1}.
 \label{eq:jp-singular}
\end{equation}
If \(f\) has bounded first and second spatial derivatives near \(x\),
then, for sufficiently small \(h\),
\begin{equation}
 \left|J_p(f(x)-f(x+h))+J_p(f(x)-f(x-h))\right|
 \le C|h|^{a_p},
 \label{eq:symmetric-smooth}
\end{equation}
where \(a_p=p\) for \(p\ge2\) and \(a_p=2(p-1)\) for \(1<p<2\).
If \(p>p_*\), its local principal value exists and the contribution
of \(|h|<\delta\) is bounded by \(C\delta^{a_p-sp}\), uniformly
on compact sets with uniform derivative bounds.
\end{lemma}

\begin{proof}
The first inequality follows by integrating \(J_p'\) along the segment
from \(a\) to \(b\). For \(1<p<2\), the concavity of
\(r^{p-1}\) on the positive axis, followed by a separation of
the two signs, gives \eqref{eq:jp-singular}.
Taylor's formula gives
\(f(x)-f(x\pm h)=\mp Df(x)\cdot h+O(|h|^2)\).
The oddness of \(J_p\), followed by \eqref{eq:jp-estimate}
or \eqref{eq:jp-singular}, gives \eqref{eq:symmetric-smooth}.
The radial integral converges when \(a_p-sp>0\), equivalently
\(p>p_*\) in the singular range.
\end{proof}

For smaller \(p\), angular cancellation is needed at a noncritical
point, and higher-order flatness is needed at a critical point.

\begin{lemma}
\label{lem:singular-pv}
Let \(1<p<2\). For a \(C^2\) spatial function with uniformly
bounded derivatives and gradient bounded away from zero, the local
principal values over \(B_\delta(x)\) are bounded by
\(C\delta^{p-sp}\), uniformly on smaller compact sets.
The same estimate holds at a point with a nonzero linear term
and a remainder bounded by \(C|h|^2\), without continuity of
the second derivatives.
If a \(C^2\) function \(f\) satisfies
\eqref{eq:critical-admissibility} at \(x_0\), its local principal
values tend to zero uniformly for centers near \(x_0\).
Consequently its patched operator, with a fixed tail in \(X\),
is continuous near \(x_0\).
\end{lemma}

\begin{proof}
The scalar estimate
\[
 |J_p(a+b)-J_p(a)|\le C_p |b|(|a|+|b|)^{p-2}
\]
follows by integrating \(J_p'\); the value at \(a=b=0\)
is understood as zero. Since
\(B\mapsto B(A+B)^{p-2}\) is increasing for \(A,B\ge0\),
a quadratic Taylor remainder can be replaced by its upper bound.
For the linear term \(Df(x)\cdot h\), spherical integration
then uses
\[
 \int_{S^{N-1}}|e\cdot\vartheta|^{p-2}\dd\vartheta<\infty
 \quad(|e|=1).
\]
The resulting radial integrand is bounded by
\(C|Df(x)|^{p-2}r^{p-sp-1}\).
This proves the noncritical assertions, including the one for a
bounded quadratic remainder. The uniform assertion near an
isolated flat critical point is
\cite[Lemmas~3.6--3.8]{KKL19}, applied on a smaller ball where
\(x_0\) is the only critical point.
The separated tail is continuous by dominated convergence and
the estimate in Lemma~\ref{lem:tail-continuity} below.
\end{proof}

\begin{lemma}
\label{lem:tail-continuity}
For \(p\ge2\), \(g,h\in X\), and \(|a|,|b|\le M\),
\begin{equation}
 \int_{\R^N}|J_p(a-g)-J_p(b-h)|\varpi\,\dd y
 \le C(1+M+\|g\|_X+\|h\|_X)^{p-2}
       \bigl(|a-b|+\|g-h\|_X\bigr),
 \label{eq:tail-nemytski}
\end{equation}
where \(C=C(N,s,p)\). For \(1<p<2\),
\begin{equation}
 \int_{\R^N}|J_p(a-g)-J_p(b-h)|\varpi\dd y
 \le C\left(|a-b|^{p-1}+d_X(g,h)\right).
 \label{eq:tail-singular}
\end{equation}
For every \(p>1\), if \(U\in C_\loc(I;X)\) and
\(J\Subset I\) is compact, then
\begin{equation}
 \sup_{\substack{t,t+\tau\in J\\ |h|+|\tau|\le r}}
 \|U(\cdot+h,t+\tau)-U(\cdot,t)\|_X\longrightarrow0
 \quad\text{as }r\downarrow0.
 \label{eq:tail-translation}
\end{equation}
\end{lemma}

\begin{proof}
For \(p>2\), apply \eqref{eq:jp-estimate} and H\"older's inequality
with exponents \((p-1)/(p-2)\) and \(p-1\) in the finite measure
\(\varpi(y)\dd y\). For \(p=2\), \eqref{eq:tail-nemytski} follows
directly from \(J_2(a)=a\).
For \(1<p<2\), \eqref{eq:tail-singular} follows from
\eqref{eq:jp-singular} and
\(|a+b|^{p-1}\le |a|^{p-1}+|b|^{p-1}\).
The translation operators \(T_hf=f(\cdot+h)\) are uniformly bounded
on \(X\) for \(|h|\le1\), since
\(\varpi(y-h)\le C\varpi(y)\). Density of \(C_c^\infty(\R^N)\)
in \(X\) gives strong continuity of \(T_h\) at \(h=0\).
This convergence is uniform on the compact set
\(\{U(\cdot,t):t\in J\}\subset X\). Combining it with uniform
time continuity on \(J\) proves \eqref{eq:tail-translation}.
When \(p<2\), these compactness and continuity statements are
made in the metric \(d_X\); they use no norm triangle inequality.
\end{proof}

\begin{remark}
For \(1<p<2\), the continuous-tail assumption can equivalently
be written in the Banach space of the evolution variable:
\[
 \beta(u)\in C_\loc\bigl(I;L^1(\R^N,\varpi\dd x)\bigr).
\]
Indeed, with \(q=p-1\), \eqref{eq:jp-singular} gives
\(\|\beta(f)-\beta(g)\|_{L^1(\varpi)}\le C d_X(f,g)\).
Conversely, the power estimate for \(\gamma\) and H\"older's
inequality give
\[
 d_X(\gamma(v),\gamma(w))
 \le C\bigl(\|v\|_{L^1(\varpi)}+\|w\|_{L^1(\varpi)}\bigr)^{1-q}
       \|v-w\|_{L^1(\varpi)}^q.
\]
Thus either description gives the same time-continuous tail class.
\end{remark}

\begin{remark}
\label{rem:tail}
The stronger assumption
\(u\in C(\R^N\times I)\cap L^\infty(\R^N\times I)\)
implies \eqref{eq:admissible-class}, by dominated convergence in
the finite measure \(\varpi(y)\dd y\). The converse does not hold:
exterior values may be unbounded or spatially discontinuous.
For every compact \(J\Subset I\), continuity into \(X\) also gives
\[
 \sup_{t\in J}\int_{|y|>R}|u(y,t)|^{p-1}\varpi(y)\dd y
 \longrightarrow0\qquad(R\to\infty).
\]
Indeed, this is uniform integrability on the compact image of \(J\)
in \(X\). A merely essentially bounded tail in time does not
specify the values on individual time slices; this distinction is
relevant to both viscosity definitions, as shown in
Section~\ref{sec:scope}.
\end{remark}

\section{Multiplicative regularization}
\label{sec:envelopes}

The homogeneity of \eqref{eq:main} suggests multiplicative
regularizations. We regularize only in a bounded interior region
and retain the original exterior values. Continuity in \(X\)
then controls the error in the equation.

The logarithms of these envelopes are inf- and sup-convolutions
of \(\log u\), defined only where \(u>0\).
A related logarithmic inf-convolution,
followed by exponentiation, was used in
\cite[proof of Theorem~8.1]{JL05}.
Separate spatial and temporal penalties ensure that a zero
spatial gradient forces the attaining spatial shift to vanish.
The regularity estimates follow from the usual competitor
argument; see \cite[Section~3 and Appendix~A]{JJ12}.

For the local Trudinger equation, the logarithmic change of variables
was established in \cite[Lemma~2.1]{BhattacharyaMarazzi15}.
Additive regularizations appear in \cite[Section~2.1]{LPS26} and,
for the fractional \(p\)-caloric equation, in
\cite[Section~3.3, Lemma~3.4]{JSU26}.

Throughout Sections~\ref{sec:envelopes}--\ref{sec:comparison},
\(u\) satisfies \eqref{eq:class}, unless stated otherwise.
At positive contact points we use the equivalent \(u\)-variable
viscosity inequalities described in Section~\ref{sec:framework}.
Fix bounded open sets \(D_0\Subset D\) and intervals
\(J_0\Subset J_1\Subset I\), and a cylinder \(Q_1\Subset Q\) such that
\(\overline D\times\overline{J_1}\subset Q_1\).
Write \(Q_0=D_0\times J_0\). Let
\(M=\max_{\overline{Q_1}}u\) and
\(m=\min_{\overline{Q_1}}u>0\). On \(\overline{Q_1}\), set
\(v=\log u\), and choose the modulus
\[
 \omega(r)=\sup\{|v(z)-v(z')|:z,z'\in\overline{Q_1},\ |z-z'|\le r\}.
\]
Here \(z=(x,t)\), and \(|\cdot|\) denotes the Euclidean norm on
\(\R^{N+1}\). Choose
\[
 \kappa=2\quad\hbox{if }p>\frac{2}{2-s},\qquad
 \kappa>\max\left\{2,\frac{sp}{p-1}\right\}
       \quad\hbox{if }1<p\le\frac{2}{2-s}.
\]
In either case \(\kappa(p-1)>sp\).
For sufficiently small \(r>0\), set
\begin{equation}
 \theta_r=2\omega(\sqrt2r)+r,\qquad
 \lambda_r=\frac{\theta_r}{r^\kappa},\qquad
 \mu_r=\frac{\theta_r}{r^2}.
 \label{eq:parameters}
\end{equation}
For \(h=(h_x,h_t)\in\R^N\times\R\), put
\[
 P_r(h)=\lambda_r|h_x|^\kappa+\mu_r h_t^2,\qquad
 \mathcal H_r=\{h:|h_x|\le r,\ |h_t|\le r\}.
\]
For \(z\in D\times J_0\), define
\begin{align}
 u_r(z)&=\min_{h\in\mathcal H_r}e^{P_r(h)}u(z+h),
 \label{eq:lower-envelope}\\
 u^r(z)&=\max_{h\in\mathcal H_r}e^{-P_r(h)}u(z+h).
 \label{eq:upper-envelope}
\end{align}
All shifted points in these formulas belong to \(Q_1\), provided
\(r\) is sufficiently small. Extend both envelopes by
\(u_r=u^r=u\) on \((\R^N\setminus D)\times J_0\).
The extensions may be discontinuous across \(\partial D\).
All regularity estimates and viscosity tests below are taken
in \(Q_0\), away from this interface.

\begin{lemma}
\label{lem:envelopes}
The extrema in \eqref{eq:lower-envelope}--\eqref{eq:upper-envelope}
are attained, and the envelopes \(u_r,u^r\)
are continuous on \(D\times J_0\).
The extensions satisfy \(u_r\le u\le u^r\) everywhere,
and \(m\le u_r,u^r\le M\) on \(D\times J_0\).
Moreover, \(u_r,u^r\in C_\loc(J_0;X)\), and
\(u_r-u\) and \(u^r-u\) converge uniformly to zero on
\(\R^N\times J_0\). In particular,
\begin{equation}
 \sup_{t\in J_0}\int_{\R^N}
 \frac{|u_r(x,t)-u(x,t)|^{p-1}+|u^r(x,t)-u(x,t)|^{p-1}}
 {1+|x|^{N+sp}}\dd x\longrightarrow0.
 \label{eq:tail-convergence}
\end{equation}
On \(Q_0\), \(u_r\) is locally Lipschitz and semiconcave, while
\(u^r\) is locally Lipschitz and semiconvex. The regularity constants
may depend on \(r\).
\end{lemma}

\begin{proof}
Compactness of \(\mathcal H_r\) and continuity of \(u\)
give attainment of the extrema and continuity
of the envelopes \(u_r,u^r\). The zero shift gives
\(u_r\le u\le u^r\), and positivity on \(Q_1\) gives
\(m\le u_r,u^r\le M\). The shifts tend uniformly to zero, while
the multipliers lie in \([1,e^{2\theta_r}]\) and
\([e^{-2\theta_r},1]\), respectively. Uniform continuity on
\(\overline{Q_1}\) and \(\theta_r\to0\) therefore give uniform
convergence on \(D\times J_0\). Since both differences vanish
outside \(D\), the convergence holds on the full strip, and
\[
 \begin{aligned}
 &\sup_{t\in J_0}\int_{\R^N}
 \frac{|u_r-u|^{p-1}+|u^r-u|^{p-1}}{1+|x|^{N+sp}}\dd x\\
 &\qquad\le \left(\int_D\varpi\,\dd x\right)
       \sup_{D\times J_0}
       \bigl(|u_r-u|^{p-1}+|u^r-u|^{p-1}\bigr).
 \end{aligned}
\]
This proves \eqref{eq:tail-convergence}. Uniform time continuity
of the interior envelopes and \(u\in C_\loc(J_0;X)\) give
\(u_r,u^r\in C_\loc(J_0;X)\), with the metric \(d_X\)
used when \(p<2\).

For \(z\in Q_0\) and \(r\) sufficiently small, all shifted points
lie in \(Q_1\), where \(u>0\). Hence
\[
 \begin{aligned}
 v_r(z):=\log u_r(z)
 &=\min_{h\in\mathcal H_r}\{v(z+h)+P_r(h)\},\\
 v^r(z):=\log u^r(z)
 &=\max_{h\in\mathcal H_r}\{v(z+h)-P_r(h)\}.
 \end{aligned}
\]
If \(h_z\) is a minimizer, comparison with the zero shift gives
\begin{equation}
 \begin{gathered}
 P_r(h_z)\le v(z)-v(z+h_z)\le\omega(\sqrt2r)<\frac{\theta_r}{2},\\
 |(h_z)_x|<2^{-1/\kappa}r,\qquad
 |(h_z)_t|<2^{-1/2}r.
 \end{gathered}
 \label{eq:inactive}
\end{equation}
The same bounds hold for maximizing shifts. For all sufficiently
small increments \(k\), the point \(z+h_z\) remains an admissible
competitor at \(z+k\) and \(z-k\). The Hessian of \(P_r\) on the
parameter region is bounded by
\(C_\kappa(\lambda_r r^{\kappa-2}+\mu_r)
=C_\kappa\theta_r/r^2\). The competitor inequalities give
\begin{equation}
 \begin{aligned}
 v_r(z+k)+v_r(z-k)-2v_r(z)
     &\le C_\kappa\frac{\theta_r}{r^2}|k|^2,\\
 v^r(z+k)+v^r(z-k)-2v^r(z)
     &\ge-C_\kappa\frac{\theta_r}{r^2}|k|^2.
 \end{aligned}
 \label{eq:semiconcave-log}
\end{equation}
Comparing neighboring points in both directions, and using the
gradient bound for \(P_r\), also gives
\begin{equation}
 |v_r(z+k)-v_r(z)|+|v^r(z+k)-v^r(z)|
 \le C_\kappa\frac{\theta_r}{r}|k|.
 \label{eq:envelope-lipschitz}
\end{equation}
Here the centers lie in the interior cylinder, and \(|k|\) is
smaller than a fixed multiple of the margins in \eqref{eq:inactive}.
Taylor's formula for the exponential, together with
\eqref{eq:semiconcave-log} and \eqref{eq:envelope-lipschitz}, gives
\[
 \begin{split}
 u_r(z+k)+u_r(z-k)-2u_r(z)&\le C_r|k|^2,\\
 u^r(z+k)+u^r(z-k)-2u^r(z)&\ge-C_r|k|^2,\\
 C_r&=C_\kappa M\frac{\theta_r+\theta_r^2}{r^2}.
 \end{split}
\]
The second inequality also follows from convexity of the exponential.
\end{proof}

The envelopes inherit the viscosity inequalities up to an error
caused by retaining the unshifted values outside \(D\).

\begin{lemma}
\label{lem:preserve}
There are numbers \(\varepsilon_r\ge0\), with
\(\varepsilon_r\to0\) as \(r\downarrow0\), such that the following holds.
If \(u\) is a viscosity supersolution in \(Q\), then every admissible
lower \(u\)-variable test for \(u_r\) in \(Q_0\) satisfies
\begin{equation}
 \beta'(u_r(z_0))\phi_t(z_0)
       +\Lp\phi_{\rho,u_r}(z_0)\ge-\varepsilon_r.
 \label{eq:approx-viscosity}
\end{equation}
For a viscosity subsolution \(u\), every admissible upper test for \(u^r\)
satisfies the reverse inequality with right-hand side
\(\varepsilon_r\). The error is uniform over the contact points
and test functions.
\end{lemma}

\begin{proof}
Suppose first that \(u\) is a viscosity supersolution.
For \(a>0\) and a space-time shift \(h\), the function
\(g(z)=a u(z+h)\) is a viscosity supersolution on the shifted
domain. Indeed, translate the test function and divide it by \(a\),
using
\begin{equation}
 \partial_t\beta(au)+\Lp(au)
 =a^{p-1}\bigl(\partial_t\beta(u)+\Lp u\bigr).
 \label{eq:homogeneity}
\end{equation}
The kernel in \eqref{eq:operator} is translation invariant.

At a contact point \(z_0\) for \(u_r\), fix an attaining shift
\(h_0\) and set \(g(z)=e^{P_r(h_0)}u(z+h_0)\).
Then \(u_r\le g\) on \(D\times J_0\) and
\(u_r(z_0)=g(z_0)\). A test \(\phi\) from below for \(u_r\) is
therefore a test from below for \(g\) on the same backward cylinder.
After translation and division by the positive multiplier, this is a
test for \(u\) at \(z_0+h_0\). For small \(r\), the translated
contact cylinder lies in \(Q\).
The temporal coefficients for \(u_r\) and \(g\) at \(z_0\) are
identical. Their patched functions also agree inside the contact ball;
on \(D\) outside that ball, the ordering
\(u_r(y,t_0)\le g(y,t_0)\) holds. Thus, with
\(c=u_r(z_0)=g(z_0)\), monotonicity of \(J_p\) gives
\[
 \Lp\phi_{\rho,u_r}(z_0)\ge\Lp\phi_{\rho,g}(z_0)
 -2\int_{\R^N\setminus D}
 \frac{|J_p(c-u(y,t_0))-J_p(c-g(y,t_0))|}
 {|x_0-y|^{N+sp}}\dd y.
\]
For \(x_0\in D_0\) and \(y\notin D\), the kernel is bounded by
\(C\varpi(y)\). Also \(|c|\le M\), and the weighted tail moments
of \(u(t_0)\) and \(g(t_0)\) are uniformly bounded. By
Lemma~\ref{lem:tail-continuity}, the last integral tends to zero
uniformly. More precisely, set
\[
 \Xi_r=\sup_{\substack{t\in J_0,\ h\in\mathcal H_r\\
                 e^{-2\theta_r}\le a\le e^{2\theta_r}}}
 \int_{\R^N}|u(y,t)-a\,u(y+h_x,t+h_t)|^{p-1}
       \varpi(y)\dd y.
\]
The shifted times lie in \(J_1\). Equation~\eqref{eq:tail-translation},
\(\theta_r\to0\), and the uniform tail bound imply \(\Xi_r\to0\).
The error is at most \(C\Xi_r\) for \(1<p<2\), by
\eqref{eq:tail-singular}, and at most \(C\Xi_r^{1/(p-1)}\) for
\(p\ge2\). We choose \(\varepsilon_r\) to bound this error for both
envelopes.
This proves \eqref{eq:approx-viscosity}. For the upper envelope,
use a maximizing shift, reverse the order on \(D\), and estimate
the same exterior error.
\end{proof}

\begin{remark}
\label{rem:amplitude}
The choice \eqref{eq:parameters} gives both the interior optimizer
estimate \eqref{eq:inactive} and \(\theta_r\to0\), which makes the
amplitude error in the exterior integral vanish. A full-space
supremal envelope need not belong to \(X\), even when the original
function does; see Section~\ref{sec:scope}.
\end{remark}

\section{Viscosity solutions are weak solutions}
\label{sec:regularized}

We first derive weak inequalities for the envelopes and an energy
bound independent of their derivative constants. We then pass to
the limit in the weak formulation.

\subsection{Weak inequalities for the envelopes}

\begin{lemma}
\label{lem:semiweak}
Let \(p>1\), and let the envelopes be defined by
\eqref{eq:lower-envelope}--\eqref{eq:upper-envelope}.
If \(u\) is a viscosity supersolution satisfying \eqref{eq:class},
then \(w=u_r\) satisfies
\[
 -\iint_{Q_0}\beta(w)\eta_t\dd x\dd t
 +\int_{J_0}\E(w,\eta)\dd t
 \ge-\varepsilon_r\iint_{Q_0}\eta\dd x\dd t
\]
for every nonnegative \(\eta\in C_c^\infty(Q_0)\).
For a viscosity subsolution \(u\), the upper envelope \(u^r\)
satisfies the reverse inequality with right-hand side
\(\varepsilon_r\iint_{Q_0}\eta\).
\end{lemma}

\begin{proof}
Fix \(r\) and let \(w=u_r\). Local Lipschitz regularity gives
\(w\in L^p_\loc(J_0;W^{s,p}_\loc(D_0))\). We begin with
a lower bound for the truncated operators.

Fix a compact subcylinder of \(Q_0\). At \((x,t)\), choose an
attaining point \((y,\tau)\) in the envelope formula and put
\(a=x-y\), \(c=w(x,t)\). For small spatial increments \(h\), the
same attaining point is an admissible competitor at \((x+h,t)\),
by \eqref{eq:inactive}. Hence
\begin{equation}
 w(x+h,t)\le F_{a,c}(h)
 :=c\exp\{\lambda_r(|a+h|^\kappa-|a|^\kappa)\}.
 \label{eq:envelope-support}
\end{equation}
The allowed size of \(h\) is uniform on the compact cylinder.
There is a constant \(C_r\) such that, for all sufficiently small
\(R>0\),
\begin{equation}
 \sup_{a,c}\sup_{0<\delta<R}
 \left|\int_{\delta<|h|<R}
 \frac{J_p(c-F_{a,c}(h))}{|h|^{N+sp}}\dd h\right|
 \le C_r\bigl(R^{p-sp}+R^{\kappa(p-1)-sp}\bigr).
 \label{eq:support-integral}
\end{equation}
Here \(|a|\le2^{-1/\kappa}r\) and \(m\le c\le M\).
For \(p\ge2\), this follows from Lemma~\ref{lem:operator}, since
the profiles \(F_{a,c}\) have uniform first and second derivative
bounds.

For \(1<p<2\), we give the estimate near \(a=0\) explicitly.
If \(a\ne0\), then
\[
 |\nabla F_{a,c}(0)|\asymp_r |a|^{\kappa-1},\qquad
 \sup_{|h|\le c_0|a|}|D^2F_{a,c}(h)|
       \le C_r|a|^{\kappa-2},
\]
where \(c_0\in(0,1/2)\) can be chosen uniformly.
The scalar inequality
\[
 |J_p(A+B)-J_p(A)|
 \le C_p|B|(|A|+|B|)^{p-2}
\]
and Taylor's formula reduce the angular integral of the symmetric
increments to an integral of
\((|e\cdot\nu|+b)^{p-2}\) over the unit sphere.
This integral is bounded uniformly for unit vectors \(\nu\) and
\(0\le b\le1\), since \(p-2>-1\). Consequently, for
\(0<\rho<c_0|a|\),
\[
 \int_{\mathbb S^{N-1}}
 \left|J_p(c-F_{a,c}(\rho e))
       +J_p(c-F_{a,c}(-\rho e))\right|\dd e
 \le C_r|a|^{\kappa(p-1)-p}\rho^p.
\]
For \(\rho\ge c_0|a|\), the bound
\(|F_{a,c}(\pm\rho e)-c|\le C_r\rho^\kappa\) applies instead.
It also applies for every \(\rho\) when \(a=0\).
To integrate the first estimate uniformly in \(a\ne0\), set
\(\alpha=\kappa(p-1)-p\).
Its radial contribution is bounded by
\[
 C_r|a|^{\alpha}
 \min\{R,c_0|a|\}^{p-sp}.
\]
If \(|a|\le R/c_0\), this is at most
\(C_r|a|^{\kappa(p-1)-sp}
\le C_rR^{\kappa(p-1)-sp}\).
If \(|a|>R/c_0\), then
\[
 |a|^{\alpha}R^{p-sp}\le
 \begin{cases}
 C R^{\kappa(p-1)-sp},
     &\alpha<0,\\
 C_rR^{p-sp},
     &\alpha\ge0,
 \end{cases}
\]
where the second case uses the fixed bound \(|a|\le r\).
The remaining radial integral is bounded by
\(C_rR^{\kappa(p-1)-sp}\), since \(\kappa(p-1)>sp\).
The same bound applies when \(a=0\), proving
\eqref{eq:support-integral}.

By monotonicity, \eqref{eq:envelope-support} and
\eqref{eq:support-integral} bound the local contribution to
\[
 \Lp^\delta w(x,t)=2\int_{|x-y|>\delta}
 \frac{J_p(w(x,t)-w(y,t))}{|x-y|^{N+sp}}\dd y
\]
from below, uniformly in \(\delta\). The remaining contribution
is bounded below by local boundedness and the uniform tail bound.
Thus \(\Lp^\delta w\ge-C_r\) on the compact cylinder.

By Alexandrov's theorem, \(w\) has a second-order
expansion at almost every space-time point; see \cite{CS04}.
Fix such a point \(z_0=(x_0,t_0)\). Suppose first that
\(D_xw(z_0)\ne0\), and let \(P\) be the Taylor polynomial at \(z_0\).
For each \(\nu>0\), the polynomial
\(P(z)-\nu|z-z_0|^2\) touches \(w\) from below in a sufficiently
small cylinder. This is an admissible test
with the same first derivatives as \(w\) at \(z_0\).
The principal value of \(w\) exists at this point.
For \(p\ge2\), this follows from the estimate
in the proof of Lemma~\ref{lem:operator}, applied to the
second-order expansion at \(z_0\).
For \(1<p<2\), the angular cancellation used above, now with
a fixed nonzero linear term and an \(O(|h|^2)\) remainder,
gives an integrable bound \(C_{z_0}\rho^p\) for the
symmetric angular integral.
The same estimates show that the difference between the operator
on the patched polynomial and the operator on \(w\) tends to zero
as the patching radius decreases. Lemma~\ref{lem:preserve}
therefore gives the desired inequality at \(z_0\).

Suppose next that \(D_xw(z_0)=0\). Differentiating the smooth
upper support supplied by an attaining point in the logarithmic
formula gives
\[
 D_x\log w(z_0)=\lambda_r\kappa
       |x_0-y|^{\kappa-2}(x_0-y),\qquad
 \partial_t\log w(z_0)=2\mu_r(t_0-\tau).
\]
Thus \(y=x_0\). Put \(b=\tau-t_0\), \(c=w(z_0)\), and
\[
 g(\xi,\sigma)=e^{\mu_r b^2}u(\xi,\sigma+b).
\]
For \((\xi,\sigma)\) near \(z_0\), the shift
\((\xi-x_0,\sigma-t_0+b)\) is admissible in the minimum at \(z_0\).
Consequently,
\[
 c\le
 e^{\lambda_r|\xi-x_0|^\kappa+\mu_r(\sigma-t_0+b)^2}
 u(\xi,\sigma+b).
\]
Multiplying by the appropriate exponential factor gives
\[
 \Phi(\xi,\sigma):=c\exp\!\left\{
 -\lambda_r|\xi-x_0|^\kappa
 -2\mu_r b(\sigma-t_0)-\mu_r(\sigma-t_0)^2
 \right\}\le g(\xi,\sigma).
\]
Equality holds at \(z_0\), so \(\Phi\) is a lower test for
\(g\) on a small backward cylinder, with
\(\Phi_t(z_0)=-2\mu_rbc=w_t(z_0)\).
When \(p\le2/(2-s)\), its spatial critical point is isolated and
its spatial profile belongs to \(C^2_\kappa\); hence it is an
admissible critical-point test. For larger \(p\), it is an ordinary
admissible smooth test.

Apply the viscosity supersolution inequality to \(g\). Although
\(\Phi\) need not touch \(w\), the ordering \(w\le g\) on \(D\)
and the exterior estimate in Lemma~\ref{lem:preserve} give
\[
 \beta'(c)w_t(z_0)+\Lp\Phi_{\rho,w}(z_0)\ge-\varepsilon_r.
\]
The local contribution of \(\Phi\) is bounded by
\(C_r\rho^{\kappa(p-1)-sp}\). Moreover,
\eqref{eq:envelope-support} with \(a=0\) gives
\(w(x_0+h,t_0)\le c e^{\lambda_r|h|^\kappa}\).
The negative part of the operator integrand is therefore
integrable near \(h=0\). Its improper integral exists in
\((-\infty,+\infty]\). Letting \(\rho\downarrow0\) proves, at
these points as well,
\begin{equation}
 \beta'(w)w_t+\Lp w\ge-\varepsilon_r
       \quad\hbox{almost everywhere}.
 \label{eq:ae}
\end{equation}

For \(\eta\in C_c^\infty(Q_0)\), \(\eta\ge0\), the lower bound
for the truncated operators permits Fatou's lemma and symmetric
integration:
\begin{align}
 \iint\eta\Lp w
 &\le\liminf_{\delta\downarrow0}\iint\eta\Lp^\delta w
 \nonumber\\
 &=\int_{J_0}\E(w,\eta)\dd t.
 \label{eq:fatou}
\end{align}
The testing form is absolutely integrable: near the test support,
the two Lipschitz increments give the power
\(|x-y|^{p-N-sp}\), and the separated interactions are controlled
by the tail bound. This proves the equality in \eqref{eq:fatou}.
Taking \(\eta=1\) on a prescribed compact cylinder shows that
\((\Lp w)_+\) is locally integrable. Since \((\Lp w)_-\)
is locally bounded, \(\Lp w\in L^1_\loc(Q_0)\); in particular,
the improper integral above is finite almost everywhere.
Integrating \eqref{eq:ae}, using
\(\beta'(w)w_t=\partial_t\beta(w)\), and applying
\eqref{eq:fatou} gives the weak supersolution inequality.

For \(w=u^r\), the lower support has profile
\(c\exp\{-\lambda_r(|a+h|^\kappa-|a|^\kappa)\}\).
The estimates above give a uniform upper bound for
\(\Lp^\delta w\). At a zero spatial gradient, the attaining
spatial point is again \(x_0\), and the comparator
\(g(\xi,\sigma)=e^{-\mu_r b^2}u(\xi,\sigma+b)\) has the upper test
\[
 c\exp\!\left\{
 \lambda_r|\xi-x_0|^\kappa
 +2\mu_r b(\sigma-t_0)+\mu_r(\sigma-t_0)^2
 \right\}.
\]
Its time derivative equals \(w_t(z_0)\). The upper-test argument
gives the reversed version of \eqref{eq:ae}; reverse Fatou and
the same integration by parts prove the subsolution assertion.
\end{proof}

\subsection{Uniform energy estimates}

The energy estimate uses the following cutoff inequality.

\begin{lemma}
\label{lem:cutoff}
For \(p>1\), \(a,b\ge0\), and \(\alpha,\gamma\in[0,1]\),
\[
 J_p(a-b)(a\alpha^p-b\gamma^p)
 \ge c_p|a\alpha-b\gamma|^p
       -C_p(a^p+b^p)|\alpha-\gamma|^p.
\]
\end{lemma}

\begin{proof}
Suppose first \(\alpha\ge\gamma\) and write the left side as
\(|a-b|^p\alpha^p+J_p(a-b)b(\alpha^p-\gamma^p)\).
Use \(|\alpha^p-\gamma^p|\le p\alpha^{p-1}|\alpha-\gamma|\)
and Young's inequality to absorb half the first term. Finally use
\(|a\alpha-b\gamma|\le|a-b|\alpha+b|\alpha-\gamma|\).
Interchanging the pairs treats \(\gamma>\alpha\).
\end{proof}

\begin{proposition}
\label{prop:energy}
If \(u\) is a viscosity supersolution satisfying \eqref{eq:class},
then for every \(B\times J\Subset Q_0\) there exists \(r_0>0\)
such that
\begin{equation}
 \sup_{0<r<r_0}\int_J[u_r(\cdot,t)]_{s,p;B}^p\dd t<\infty.
 \label{eq:uniform-energy}
\end{equation}
The same statement holds for \(u^r\) when \(u\) is a viscosity
subsolution. The energy bound depends on the nested cylinders,
\(N,s,p\), a local upper bound for \(u\), and its uniform \(X\)-bound
on a compact time interval. It does not depend on a positive lower
bound or on the derivative constants of the envelopes.
The admissible threshold \(r_0\), however, may depend on the positive
lower bound and the modulus of continuity of \(\log u\).
\end{proposition}

\begin{proof}
Choose \(K=1+M\), where \(M\) is the local upper bound in
Section~\ref{sec:envelopes}, and decrease \(r_0\) so that
\(\varepsilon_r\le1\) for \(r<r_0\).
By Lemmas~\ref{lem:envelopes}, \ref{lem:preserve}, and
\ref{lem:semiweak}, \(w=u_r\) satisfies the weak supersolution
inequality with right-hand side \(-\varepsilon_r\).
Choose a cylinder \(A\times J_*\Subset Q_0\) containing
\(\overline B\times\overline J\), and cutoffs
\(0\le\zeta\le1\), \(\zeta\in C_c^\infty(A)\), and
\(0\le\chi\in C_c^\infty(J_*)\), equal to one on \(B\) and
\(J\), respectively. On \(A\times J_*\), \(0<w\le K\).
The envelope bounds and equality with \(u\) outside \(D\) give
\[
 \mathcal T:=\sup_{0<r<r_0}\sup_{t\in J_*}
       \bigl(\|u_r(\cdot,t)\|_X+\|u^r(\cdot,t)\|_X\bigr)<\infty.
\]
Test with \((K-w)\zeta^p\chi\). For fixed \(r\),
\(\partial_t\beta(w)=\beta'(w)w_t\) is locally bounded, and the test is nonnegative with compact support in the
local energy space. Lemma~\ref{lem:energy-testing} applies, with
the constant source term handled by the same approximation. For \(a>0\), set
\[
 H(a)=K a^{p-1}-\frac{p-1}{p}a^p,
 \qquad H'(a)=\beta'(a)(K-a).
\]
The time integral becomes
\(-\iint H(w)\zeta^p\chi'\), and \(|H(w)|\le C_pK^p\).
Put \(V=K-w\) only on \(A\), and extend \(Z=V\zeta\) by zero
outside \(A\). Split the testing form into \(A\times A\) and
the pairs joining \(\supp\zeta\) to \(A^c\).
On \(A\times A\), \(0\le V\le K\), so
Lemma~\ref{lem:cutoff} gives
\[
 \iint_{A\times A}
 \frac{J_p(w(x)-w(y))(V(x)\zeta(x)^p-V(y)\zeta(y)^p)}
 {|x-y|^{N+sp}}\dd x\dd y
 \le-c_p[Z]_{s,p;A}^p+C_pK^p[\zeta]_{s,p}^p.
\]
In the remaining pairs, the positive separation and the tail bound
give an absolute bound
\[
 C K(K^{p-1}+\mathcal T^{p-1})\|\zeta\|_{L^p}^p.
\]
The difference between \([Z]_{s,p}^p\) and \([Z]_{s,p;A}^p\)
is bounded by \(CK^p\|\zeta\|_{L^p}^p\), again by separation.
Finally, the source term contributes at most
\(\varepsilon_r K\|\zeta\|_{L^p}^p\int\chi\).
Combining these estimates yields
\begin{equation}
 \begin{split}
 \int\chi[Z]_{s,p}^p\dd t
 \le{}&C K^p\left(
 \|\zeta\|_{L^p}^p\int|\chi'|\dd t
 +[\zeta]_{s,p}^p\int\chi\dd t\right)\\
 &+C\bigl(K^p+K\mathcal T^{p-1}+\varepsilon_r K\bigr)
        \|\zeta\|_{L^p}^p\int\chi\dd t .
 \end{split}
 \label{eq:caccioppoli}
\end{equation}
Taking cutoffs equal to one on \(B\times J\) proves
\eqref{eq:uniform-energy}.

The constant \(C\) depends only on \(N,s,p\) and the fixed
spatial geometry and separation. Since \(\varepsilon_r\le1\), it is independent
of \(r\) and of the positive lower bound.
For a subsolution, use \((w+K)\zeta^p\chi\).
Now \(V=w+K\) lies in \([0,2K]\) on \(A\), the spatial form has
the positive sign, and the time primitive is
\(K a^{p-1}+(p-1)a^p/p\). The same splitting and source estimate
give the asserted bound.
\end{proof}

\subsection{Passage to the limit}
\label{sec:limit}

The uniform energy bound gives local Sobolev regularity of the
limit. The local flux converges weakly, while the separated
interactions converge by dominated convergence.

\begin{proposition}
\label{prop:viscweak}
A viscosity supersolution satisfying \eqref{eq:class} is a local weak
supersolution. The analogous assertion holds for subsolutions.
\end{proposition}

\begin{proof}
Take a sequence \(r_j\downarrow0\) and use the lower envelopes in
an arbitrary system of nested cylinders as in
Section~\ref{sec:envelopes}. Local uniform
convergence and Proposition~\ref{prop:energy} imply
\(u\in L^p(J;W^{s,p}(A))\) for every
\(A\times J\Subset Q_0\), by Fatou's lemma.

To pass to the testing form, choose \(D_*\Subset D_0\) containing the spatial
support of a test \(\eta\) at positive distance from its boundary.
On \(D_*\times D_*\times J\), use the measure
\[
 \dd\mu=|x-y|^{-N-sp}\dd x\dd y\dd t.
\]
The functions
\(G_j(x,y,t)=u_{r_j}(x,t)-u_{r_j}(y,t)\) are bounded in \(L^p(\mu)\)
and converge almost everywhere to
\(G(x,y,t)=u(x,t)-u(y,t)\). Hence \(J_p(G_j)\) is bounded in
\(L^{p'}(\mu)\). Weak compactness and almost-everywhere convergence give
\begin{equation}
 J_p(G_j)\rightharpoonup J_p(G)\quad\hbox{in }L^{p'}(\mu).
 \label{eq:flux}
\end{equation}
Indeed, on finite-measure sets separated from the diagonal,
almost-everywhere convergence and uniform integrability identify
every weak subsequential limit with \(J_p(G)\). Exhausting the
measure space by such sets and using the uniform \(L^{p'}(\mu)\)
bound proves \eqref{eq:flux}.

Pair \eqref{eq:flux} with \(\eta(x,t)-\eta(y,t)\in L^p(\mu)\).
The remaining pairs have one point in \(\supp\eta\) and the other
outside \(D_*\). Their distance is bounded below. The envelopes
are uniformly bounded on \(D\) and equal \(u\) on \(D^c\).
Thus an integrable majorant is furnished by the local bound and
\((1+|u(y,t)|^{p-1})\varpi(y)\) in the exterior.
Dominated convergence treats this part.
Finally \(\beta(u_{r_j})\to\beta(u)\) uniformly on the test support.
The right-hand sides
\(-\varepsilon_{r_j}\iint\eta\) tend to zero.
Passing to the limit in the approximate weak supersolution inequalities proves
the assertion. The upper envelopes yield the subsolution assertion.
\end{proof}

\section{Weak solutions are viscosity solutions}
\label{sec:comparison}

We prove the converse by comparing a weak solution with smooth
subsolutions and supersolutions. The comparison argument uses a
relative-energy identity for the nonlinear time derivative and
requires the weak equation with equality.

\subsection{Stability under smooth perturbations}

We first estimate the effect of perturbations supported inside the
patching ball.

\begin{lemma}
\label{lem:strict-stability}
Assume \(p>2/(2-s)\). Let \(D_0\Subset D_1\) be bounded open sets and let \(J\Subset I\)
be a compact time interval. Let \(w\) be measurable and satisfy
\(\sup_{t\in J}\|w(\cdot,t)\|_X<\infty\).
Assume that \(w\) is \(C^{2,1}\) on a neighborhood of
\(\overline D_1\times J\) and that \(w\ge m_0>0\) there.
Let \(\psi\) be smooth with compact spatial support in \(D_1\)
and bounded derivatives of the indicated orders. For small
\(|\varepsilon|\), set \(w_\varepsilon=w+\varepsilon\psi\).
Then, uniformly on \(D_0\times J\),
\begin{align*}
 |\Lp w_\varepsilon-\Lp w|
   &\le C|\varepsilon|^{\sigma_p},\\
 |\mathcal R(w_\varepsilon)-\mathcal R(w)|
   &\le C\bigl(|\varepsilon|+|\varepsilon|^{\sigma_p}\bigr),
\end{align*}
where \(\mathcal R(w)=(p-1)w^{p-2}w_t+\Lp w\) and
\[
 \sigma_p=
 \begin{cases}
  1-s,&p\ge2,\\
  p-1-sp/2,&2/(2-s)<p<2.
 \end{cases}
\]
The constant depends on \(N,s,p\), \(m_0\), the two domains,
the local smoothness and tail bounds for \(w\), and the indicated
derivative bounds for \(\psi\). It is independent of
\(\varepsilon\).
In particular, if
\(\mathcal R(w)\ge\delta\) or
\(\mathcal R(w)\le-\delta\) on \(D_0\times J\)
for some \(\delta>0\), then the corresponding strict sign
is preserved for sufficiently small \(|\varepsilon|\).
\end{lemma}

\begin{proof}
Choose \(0<r_0<\min\{1,\operatorname{dist}
(\overline D_0,\R^N\setminus D_1)/2\}\). Suppose first that
\(p\ge2\). The symmetric spatial
increments satisfy
\[
 \bigl|J_p(w(x,t)-w(x+h,t))
       +J_p(w(x,t)-w(x-h,t))\bigr|\le C|h|^p,
\]
and likewise for \(w_\varepsilon\). This follows from the
cancellation of the linear Taylor terms and
\eqref{eq:jp-estimate}. The difference of the near-field integrals
over \(|h|<\rho\) is therefore bounded by \(C\rho^{p-sp}\).
For \(\rho\le |h|<r_0\), local boundedness and
\eqref{eq:jp-estimate} give
\[
 \int_{\rho\le |h|<r_0}
 \frac{|J_p(w_\varepsilon(x,t)-w_\varepsilon(x+h,t))
          -J_p(w(x,t)-w(x+h,t))|}{|h|^{N+sp}}\dd h
 \le C|\varepsilon|\rho^{-sp}.
\]
For \(|h|\ge r_0\), the same inequality gives a bound by
\[
 C|\varepsilon|\int_{|h|\ge r_0}
 \frac{1+|w(x+h,t)|^{p-2}}{|h|^{N+sp}}\dd h
 \le C|\varepsilon|.
\]
For this last estimate, the separated kernel is bounded by
\(C\varpi(x+h)\), uniformly for \(x\in D_0\), and H\"older's
inequality applies with the uniform \(X\)-bound. When \(p=2\),
the numerator is constant.
Take \(\rho=r_0|\varepsilon|^{1/p}\) for
\(0<|\varepsilon|\le1\). All three contributions are bounded by
\(C|\varepsilon|^{1-s}\).

If \(2/(2-s)<p<2\), use
\(|J_p(a)-J_p(b)|\le C_p|a-b|^{p-1}\).
The symmetric increments of a \(C^2\) function are bounded by
\(C|h|^{2(p-1)}\), so the near-field difference is at most
\(C\rho^{2(p-1)-sp}\). The perturbation of the remaining
integral is bounded by \(C|\varepsilon|^{p-1}\rho^{-sp}\).
Choosing \(\rho=r_0|\varepsilon|^{1/2}\) gives the exponent
\(\sigma_p=p-1-sp/2>0\).
In both cases the time term changes by at most
\(C|\varepsilon|\), since \(w_\varepsilon\ge m_0/2\) for small
\(|\varepsilon|\), which proves the estimates.
\end{proof}

For \(1<p\le2/(2-s)\), a cutoff constant near the contact point
preserves the admissibility of a critical test.

\begin{lemma}
\label{lem:frozen-stability}
Let \(1<p\le2/(2-s)\), \(J\Subset I\) be compact, and
\(B_R(x_0)\Subset\Omega\). Suppose that \(w\in C(J;X)\) and
\[
 w(x,t)=f(x)+b(t)\quad\hbox{on }B_R(x_0)\times J,
 \qquad f\in C^2(B_R(x_0)),\quad b\in C^1(J).
\]
Assume that \(f\) is spatially admissible at \(x_0\) in the
sense of Section~\ref{sec:framework}.
For every sufficiently small \(r>0\), there is
\(\eta\in C_c^\infty(B_{r/2}(x_0))\) with
\(0\le\eta\le1\) and \(\eta=1\) on \(B_{r/4}(x_0)\)
such that, for every measurable \(\vartheta:J\to[0,1]\) and
sufficiently small \(|\varepsilon|\),
\[
 \sup_{B_r(x_0)\times J}
 \bigl|\Lp(w+\varepsilon\eta\vartheta)-\Lp w\bigr|
 \le C|\varepsilon|^{(p-1)(1-s)}.
\]
Here \(C\) depends on \(N,s,p\), the spatial profile, and the
fixed radii, but not on \(\varepsilon\) or \(\vartheta\).
If, in addition, \(\vartheta\in C^1(J)\) and
\(0<m_0\le w\le M_0\) on \(B_r(x_0)\times J\), then
\[
 \sup_{B_r(x_0)\times J}
 \left|\partial_t\beta(w+\varepsilon\eta\vartheta)
             -\partial_t\beta(w)\right|
 \le C_0|\varepsilon|
       \bigl(\|b'\|_{L^\infty(J)}
                    +\|\vartheta'\|_{L^\infty(J)}\bigr),
\]
where \(C_0=C_0(p,m_0,M_0)\), provided
\(|\varepsilon|\le\min\{1,m_0/2\}\).
The upper bound \(M_0\) is finite by local continuity.
\end{lemma}

\begin{proof}
Shrink \(R\) if necessary. Either \(Df\) does not vanish there,
or \(x_0\) is its only zero. Choose \(r<R/4\) and the stated
plateau cutoff. On the transition annulus, \(|Df|\) has a
positive minimum. For small \(|\varepsilon|\), the perturbed
profiles \(f+\varepsilon\vartheta(t)\eta\) therefore have no
new critical points. Near \(x_0\) their spatial derivatives
coincide with those of \(f\), so they retain its admissibility.
Their symmetric principal values converge locally uniformly by
Lemma~\ref{lem:singular-pv}; see also the plateau construction in
\cite[Lemma~3.9]{KKL19}.

For \(x\in B_{r/8}(x_0)\), the perturbation is spatially
constant on a fixed neighborhood of \(x\). Its near-field
increments cancel exactly, while the separated part is bounded by
\(C|\varepsilon|^{p-1}\), using the global
\((p-1)\)-H\"older continuity of \(J_p\).
For \(x\in B_r(x_0)\setminus B_{r/8}(x_0)\), the gradients
of all the perturbed profiles are bounded away from zero on a
fixed neighborhood. The noncritical estimate in
Lemma~\ref{lem:singular-pv} bounds the symmetric integral over
\(|h|<\rho\) by \(C\rho^{p-sp}\), uniformly in \(t\) and
\(\varepsilon\). The difference of the remaining integrals is
at most \(C|\varepsilon|^{p-1}\rho^{-sp}\).
Taking \(\rho=c_r|\varepsilon|^{(p-1)/p}\), with \(c_r>0\)
fixed and small, proves the spatial estimate. For the time term,
write \(w_\varepsilon=w+\varepsilon\eta\vartheta\). Since
\(w_t=b'\),
\[
 \partial_t\beta(w_\varepsilon)-\partial_t\beta(w)
 =\bigl[\beta'(w_\varepsilon)-\beta'(w)\bigr]b'
       +\varepsilon\beta'(w_\varepsilon)\eta\vartheta'.
\]
Both \(w\) and \(w_\varepsilon\) take values in
\([m_0/2,M_0+1]\). Bounds for \(\beta'\) and
\(\beta''\) on this interval give the stated estimate.
\end{proof}

\subsection{Composition and chain rules}

We record the composition and chain rules needed to test with the
positive part of a difference.

\begin{lemma}
\label{lem:composition}
Let \(0<s<1\), \(1\le p<\infty\), and let \(J\) be a bounded
interval. Let \(G:\R^\ell\to\R\) be globally Lipschitz with \(G(0)=0\).
If \(U_j\to U\) in \(L^p(J;W^{s,p}(\R^N;\R^\ell))\), then
\(G(U_j)\to G(U)\) in \(L^p(J;W^{s,p}(\R^N))\).
\end{lemma}

\begin{proof}
We follow the composition argument in
\cite[Lemma~2.5]{BPS14}.
From every subsequence choose a further subsequence such that
\(\sum_j\|U_j-U\|_{L^p(J;W^{s,p})}<\infty\).
Write
\[
 D_sU(x,y,t)=\frac{U(x,t)-U(y,t)}{|x-y|^{N/p+s}},
 \qquad H=|D_sU|+\sum_j|D_s(U_j-U)|.
\]
Minkowski's inequality gives \(H\in L^p(\R^{2N}\times J)\).
The difference quotients of \(G(U_j)-G(U)\) converge almost
everywhere to zero and are bounded in absolute value by
\(2\operatorname{Lip}(G)H\).
Dominated convergence gives convergence in the seminorm, while
\(L^p\) convergence follows from the Lipschitz bound. Since this
argument applies to every subsequence, the full sequence converges.
\end{proof}

In the application below, \(A=\beta(v)\) for a smooth comparison
function \(v\), and \(B\) is a smooth extension of
\(\beta^{-1}\).

\begin{lemma}
\label{lem:abstract-chain}
Let \(D\) be bounded and open, \(J=(a,b)\), \(p>1\), and \(V=V(D)\).
Suppose that
\[
 q\in L^p(J;V)\cap L^\infty(D\times J)
       \cap C([a,b];L^2(D)),\qquad
 q_t\in L^{p'}(J;V^*).
\]
Let \(A\) be \(C^1\) on a neighborhood of
\(\overline D\times[a,b]\). Let \(B\in C^2(\R)\) satisfy
\(0<c\le B'\le C\) and \(\|B''\|_\infty<\infty\).
For \(\sigma\in\{-1,1\}\), define
\begin{align*}
 g_\sigma(A,q)&=\sigma\bigl[B(A+\sigma q_+)-B(A)\bigr],\\
 H_\sigma(A,q)&=\int_0^{q_+}
       \sigma\bigl[B(A+\sigma\tau)-B(A)\bigr]\dd\tau.
\end{align*}
Then \(g_\sigma(A,q)\in L^p(J;V)\), the integral of
\(H_\sigma(A,q)\) over \(D\) is absolutely continuous on \([a,b]\),
and
\begin{equation}
 \frac{\dd}{\dd t}\int_D H_\sigma(A,q)\dd x
 =\pair{q_t}{g_\sigma(A,q)}
 +\int_D\sigma A_t
       \bigl[g_\sigma(A,q)-B'(A)q_+\bigr]\dd x
 \label{eq:abstract-chain}
\end{equation}
for almost every \(t\in J\).
\end{lemma}

\begin{proof}
The map \((A,q)\mapsto g_\sigma(A,q)\) is globally Lipschitz,
and \(g_\sigma(A,0)=0\). Extend \(A\) as a \(C^1\) function with
compact spatial support, without changing it near \(\overline D\).
Lemma~\ref{lem:composition} gives the Sobolev regularity and
strong continuity in \(q\). To obtain membership in \(V\),
approximate \(q\) by functions in \(C_c^\infty(D)\).
The compositions have compact support in \(D\) and can be
approximated in \(W^{s,p}(\R^N)\) by smooth functions.
A second application of Lemma~\ref{lem:composition} yields
\(g_\sigma(A,q)\in L^p(J;V)\).

On a compact interval \(J_0\Subset J\), take the forward Steklov average
\[
 q_h(t)=\frac{1}{h}\int_t^{t+h}q(\tau)\dd\tau.
\]
For small \(h\),
\[
 q_h\to q\text{ in }L^p(J_0;V),\qquad
 (q_h)_t\to q_t\text{ in }L^{p'}(J_0;V^*).
\]
The second convergence follows from
\((q_h)_t=(q_t)_h\) in \(V^*\).
Also \(q_h\to q\) uniformly in \(L^2(D)\) on \(J_0\), and
\(\|q_h\|_\infty\le\|q\|_\infty\).
For each fixed \(h\), the difference-quotient formula gives
\(q_h\in W^{1,p}(J_0;V)\). Moreover,
\((q_h)_t=(q(t+h)-q(t))/h\), and the assumed \(L^\infty\)
bound gives \(q_h\in W^{1,\infty}(J_0;L^2(D))\).
The ordinary Hilbert-space chain rule therefore applies without
requiring an embedding of \(V\) into \(L^2(D)\). Its derivatives are
\[
 \partial_qH_\sigma=g_\sigma,
 \qquad
 \partial_AH_\sigma
   =\sigma\bigl[g_\sigma-B'(A)q_+\bigr].
\]
For each fixed \(h\), the \(L^\infty(D)\)
bound on \((q_h)_t\) ensures that the \(L^2\) inner product
agrees with the \(V^*,V\) duality pairing:
\[
 \int_D (q_h)_t\,g_\sigma(A,q_h)\dd x
 =
 \pair{(q_h)_t}{g_\sigma(A,q_h)}.
\]
Lemma~\ref{lem:composition} gives
\(g_\sigma(A,q_h)\to g_\sigma(A,q)\) strongly in
\(L^p(J_0;V)\). Hence the dual pairings converge in \(L^1(J_0)\).
The terms involving \(A_t\) converge in \(L^1\) by the Lipschitz
bounds and boundedness of \(A_t\). The energy integrals converge
uniformly in time on \(J_0\) by the uniform \(L^2\) convergence.
Passing to the integrated identity proves
\eqref{eq:abstract-chain} on \(J_0\). Its right-hand side belongs
to \(L^1(J)\); continuity at the endpoints therefore gives the
absolutely continuous identity on \([a,b]\).
\end{proof}

\subsection{The relative-energy identity}

Let \(\beta(a)=a^{p-1}\) on a positive interval
\([m_*,M_*]\), with \(0<m_*<M_*<\infty\). Define
\begin{align*}
 F_+(a,b)&=\mathbf1_{\{a>b\}}\int_b^a(r-b)\beta'(r)\dd r,\\
 F_-(a,b)&=\mathbf1_{\{a<b\}}\int_a^b(b-r)\beta'(r)\dd r.
\end{align*}
There are constants depending only on \(p,m_*,M_*\) such that
\begin{equation}
 c(a-b)_\pm^2\le F_\pm(a,b)\le C(a-b)_\pm^2,
 \label{eq:relative-coercivity}
\end{equation}
where \((a-b)_-=(b-a)_+\).

\begin{lemma}
\label{lem:chain}
Let \(D\Subset\Omega\) be bounded and open, and let
\([a,b]\Subset I\). Let \(u\) be a local weak solution,
continuous in a neighborhood of
\(\overline D\times[a,b]\Subset Q\), and let \(v\in C^{2,1}\)
there. Suppose \(m_*\le u,v\le M_*\) on
\(\overline D\times[a,b]\). If the function \((u-v)_+\) on \(D\)
is supported in a fixed compact subset of \(D\),
uniformly for \(t\in[a,b]\), then
\begin{equation}
 \pair{\partial_t(\beta(u)-\beta(v))}{(u-v)_+}
 =\frac{\mathrm d}{\mathrm dt}\int_DF_+(u,v)\dd x
   +\int_Dv_t R_+(u,v)\dd x,
 \label{eq:chain-plus}
\end{equation}
in the integrated time sense, where
\[
 R_+(u,v)=\mathbf1_{\{u>v\}}
 [\beta(u)-\beta(v)-\beta'(v)(u-v)].
\]
If instead \((v-u)_+\) has such compact support, then
\begin{equation}
 \pair{\partial_t(\beta(v)-\beta(u))}{(v-u)_+}
 =\frac{\mathrm d}{\mathrm dt}\int_DF_-(u,v)\dd x
   +\int_Dv_t R_-(u,v)\dd x,
 \label{eq:chain-minus}
\end{equation}
where
\[
 R_-(u,v)=\mathbf1_{\{v>u\}}
 [\beta'(v)(v-u)-(\beta(v)-\beta(u))].
\]
Both remainder terms satisfy \(|R_\pm(u,v)|\le C F_\pm(u,v)\).
They are nonnegative for \(p\ge2\) and nonpositive for
\(1<p<2\).
\end{lemma}

\begin{proof}
For smooth functions, the identities follow from
\[
 \partial_uF_+(u,v)=\beta'(u)(u-v)_+,
 \qquad
 \partial_vF_+(u,v)=-(\beta(u)-\beta(v))_+.
\]
For \(F_-\), the derivatives on \(u<v\) are
\(-\beta'(u)(v-u)\) and \(\beta(v)-\beta(u)\).
To obtain the identities for \(u\), we apply
Lemma~\ref{lem:abstract-chain} after localization.

Write \(z=\beta(u)\), and choose \(0\le\xi\le1\), smooth and
compactly supported in \(D\), equal to one near the fixed support
of the positive part used in the chosen identity. For
\(\eta\in V(D)\), Lemma~\ref{lem:energy-testing} and the weak
equation give
\[
 \pair{(\xi z)_t}{\eta}=-\E(u,\xi\eta)
 \quad\hbox{in distributions in time}.
\]
On a slightly larger domain \(D'\Subset\Omega\), the near-field
contribution is bounded by
\(C[u]_{W^{s,p}(D')}^{p-1}\|\eta\|_{V(D)}\). The separated
contribution is bounded by
\[
 C\left(\|u\|_{L^\infty(D'\times[a,b])}^{p-1}
       +\|u(\cdot,t)\|_X^{p-1}\right)\|\eta\|_{V(D)}.
\]
Integration in time shows that \((\xi z)_t\) belongs to
\(L^{p'}(a,b;V(D)^*)\). Multiplication by \(\xi\) is bounded
in the fractional Sobolev norm, and \(\beta\) is Lipschitz on
the relevant bounded range. These facts also give
\(\xi z\in L^p(a,b;V(D))\).

Next choose a compact positive interval whose interior contains
\([\beta(m_*),\beta(M_*)]\). Extend \(\beta^{-1}\) from this
interval to a function \(B\in C^2(\R)\) with
\(0<c\le B'\le C\) and bounded \(B''\). Such an extension is
obtained by extending its positive derivative smoothly to positive
constants outside a larger interval and integrating.
Thus \(B\) remains defined after time regularization and agrees
with \(\beta^{-1}\) on the relevant positive interval.

Set \(A=\beta(v)\) and, for the chosen sign,
\(q=\sigma\xi(z-A)\), where \(\sigma=1\) is the plus case and
\(\sigma=-1\) is the minus case. Only values of \(A\) near
\(\supp\xi\) matter. A larger spatial cutoff therefore supplies
the \(C^1\) extension of \(A\) required by
Lemma~\ref{lem:abstract-chain}, without changing \(q\) or any
nonzero integrand below.
The preceding estimates give the two evolution-space bounds in
that lemma. Boundedness and continuity of \(u\), together with the
compact support of \(\xi\), give the boundedness and the continuous
\(L^2\) representative required there.
Where \(q>0\), the cutoff is one. Consequently
\[
 \begin{array}{c|c|c|c}
 \sigma&q_+&g_\sigma(A,q)&H_\sigma(A,q)\\ \hline
 1&(\beta(u)-\beta(v))_+&(u-v)_+&F_+(u,v)\\
 -1&(\beta(v)-\beta(u))_+&(v-u)_+&F_-(u,v).
 \end{array}
\]
The entries also agree where the respective positive part vanishes.
On its support, \(B'(A)=1/\beta'(v)\) and
\(A_t=\beta'(v)v_t\). Thus the last term in
\eqref{eq:abstract-chain} is exactly
\(-\int_Dv_tR_+(u,v)\) or \(-\int_Dv_tR_-(u,v)\).
In the dual pairing, \(\xi g_\sigma=g_\sigma\), so the cutoff
disappears. This proves \eqref{eq:chain-plus} and
\eqref{eq:chain-minus}, including their endpoint interpretation.

Finally \(\beta(a)=a^{p-1}\) is convex for \(p\ge2\) and
concave for \(1<p<2\), which gives the stated signs.
Taylor's formula gives the bound
\begin{equation}
 |R_\pm(u,v)|\le C_\beta F_\pm(u,v),\qquad
 C_\beta=
 \frac{\max_{[m_*,M_*]}|\beta''|}
      {\min_{[m_*,M_*]}\beta'}.
 \label{eq:explicit-relative-remainder}
\end{equation}
Indeed \(F_+\ge\tfrac12\min\beta'\,(u-v)_+^2\) and
\(F_-\ge\tfrac12\min\beta'\,(v-u)_+^2\), whereas each remainder
is at most \(\tfrac12\max|\beta''|\) times the corresponding
square. For \(p=2\) the remainders vanish and the energies reduce
to one half of the squared positive parts.
\end{proof}

\subsection{Comparison with smooth subsolutions and supersolutions}

\begin{lemma}
\label{lem:comparison}
Let \(D\Subset\Omega\) be bounded and open, and let
\([a,b]\Subset I\). Let \(u\) be a local weak solution,
continuous near
\(\overline D\times[a,b]\Subset Q\). Let \(v\) be measurable and
belong to \(L^\infty((a,b);X)\), and be \(C^{2,1}\) near this closed
cylinder. Assume \(0<m_*\le u,v\le M_*<\infty\) inside the cylinder.
Let \(K\Subset D\) be a fixed compact set.

Suppose that
\[
 \begin{gathered}
 \partial_t\beta(v)+\Lp v\ge0\quad\hbox{weakly in }D\times(a,b),\\
 u\le v\quad\hbox{on }(\R^N\setminus K)\times[a,b],
 \qquad u(\cdot,a)\le v(\cdot,a)\quad\hbox{in }D.
 \end{gathered}
\]
Then \(u\le v\) in \(D\times[a,b]\).
Reversing all three inequalities gives \(v\le u\) instead.
The exterior order is understood almost everywhere in space and
time; the initial order is imposed on the continuous interior
traces.
\end{lemma}

\begin{proof}
Local smoothness and the uniform \(X\)-bound make the testing
form of \(v\) finite. For the pointwise comparison functions used
below, locally uniform convergence of the symmetric principal
values also gives the weak inequality by symmetric integration.

For the first assertion set \(w=(u-v)_+\) on the whole space.
The exterior order makes it zero almost everywhere outside \(K\).
Interior continuity gives \(u\le v\) on
\((D\setminus K)\times[a,b]\) at every point, so the interior
support of \(w\) lies in \(K\) at every time. Composition with
the Lipschitz positive-part map and multiplication by a cutoff
equal to one near \(K\) give \(w\in L^p(a,b;V(D))\).
Lemma~\ref{lem:energy-testing} permits subtraction of the two
relations and pairing with \(w\) in the dual energy space.
Indeed, \(\partial_t\beta(v)\) is bounded on the cylinder and
therefore belongs to \(L^{p'}(a,b;V(D)^*)\). Both time derivatives
are thus admissible in the dual pairing.
Apply Lemma~\ref{lem:chain} on compact time intervals and pass
to the endpoints. We obtain
\[
 \frac{\mathrm d}{\mathrm dt}\int_DF_+(u,v)\dd x
 +\E(u,w)-\E(v,w)
 \le-\int_Dv_tR_+(u,v)\dd x.
\]
For each pair \((x,y)\), let \(d(x)=u(x)-v(x)\). The sign of
\[
 J_p(u(x)-u(y))-J_p(v(x)-v(y))
\]
is the sign of \(d(x)-d(y)\). The same ordering holds for
\(d_+(x)-d_+(y)\), or that difference is zero. Hence the integrand
in \(\E(u,w)-\E(v,w)\) is nonnegative. All these testing integrals
are finite by local energy and the uniform weighted tail bounds.

Put \(Y_+(t)=\int_DF_+(u,v)\dd x\).
Estimate~\eqref{eq:explicit-relative-remainder} gives
\[
 Y_+'(t)
 \le C_\beta\|v_t(\cdot,t)\|_{L^\infty(D)}Y_+(t)
 \quad\hbox{for almost every }t.
\]
The energy is absolutely continuous by Lemma~\ref{lem:chain}.
Gronwall's inequality gives
\[
 Y_+(t)\le Y_+(a)
 \exp\left(C_\beta\int_a^t
       \|v_t(\cdot,\tau)\|_{L^\infty(D)}\dd\tau\right).
\]
The initial order makes \(Y_+(a)=0\), and
\eqref{eq:relative-coercivity} gives \(w=0\). Interior continuity
extends the order to every point and to the terminal time.
For the reversed inequalities use \(w=(v-u)_+\) and
\eqref{eq:chain-minus}; the spatial difference is then
\(\E(v,w)-\E(u,w)\ge0\). The same estimate holds with
\(Y_-(t)=\int_DF_-(u,v)\dd x\), and proves the second assertion.
\end{proof}

\subsection{Proof of the weak-to-viscosity implication}

The comparison cylinder will be smaller than the contact cylinder.
We keep the original patch outside it, so that the exterior values
in the viscosity test remain unchanged.

\begin{proposition}
\label{prop:weakvisc}
A weak solution satisfying \eqref{eq:class} is a viscosity solution
in \(Q\).
\end{proposition}

\begin{proof}
For \(1<p<2\), Lemma~\ref{lem:positive-tests} reduces the proof
to tests for \(u\) satisfying the spatial admissibility condition
at the contact time. Let \(\phi\) be a test from below at
\(z_0=(x_0,t_0)\), with patching radius \(\rho\). Suppose that
\[
 \beta'(u(z_0))\phi_t(z_0)+\Lp\phi_{\rho,u}(z_0)=-4\delta<0.
\]
Put \(a=\phi_t(z_0)\) and choose \(\nu>0\) with
\(\beta'(u(z_0))\nu<\delta\). By continuity of \(\phi_t\),
there are \(R<\rho\) and \(\tau>0\) such that
\(\lvert\phi_t-a\rvert<\nu\) on
\(B_R(x_0)\times[t_0-\tau,t_0]\). Define
\[
 \Phi(y,t)=
 \begin{cases}
  \phi(y,t_0)+(a+\nu)(t-t_0),&y\in B_R(x_0),\\
  \phi_{\rho,u}(y,t),&y\notin B_R(x_0).
 \end{cases}
\]
For backward times, integration of \(\phi_t\le a+\nu\)
gives \(\Phi\le\phi_{\rho,u}\le u\). At time \(t_0\) the
entire spatial patch is unchanged. Consequently
\[
 \partial_t\beta(\Phi)(z_0)+\Lp\Phi(z_0)<-3\delta.
\]
Inside \(B_R(x_0)\), the spatial derivatives of \(\Phi\) are
independent of time. The admissibility of a critical contact profile
is therefore preserved at nearby times. Lemmas~\ref{lem:operator} and
\ref{lem:singular-pv}, together with tail continuity, show that
the last residual is continuous near \(z_0\). Hence there are
\(t_*<t_0\) and a ball \(D_1\Subset B_R(x_0)\) on which
\(u,\Phi\) lie in a compact positive interval and
\[
 \partial_t\beta(\Phi)+\Lp\Phi\le-2\delta
       \quad\hbox{on }\overline D_1\times[t_*,t_0].
\]

Choose a smaller ball \(D=B_r(x_0)\Subset D_1\) and a cutoff
\(\eta\in C_c^\infty(D)\), with \(0\le\eta\le1\) and
\(\eta(x_0)=1\). If \(p\le2/(2-s)\), choose \(D\) and
the plateau cutoff as in Lemma~\ref{lem:frozen-stability}.
Let \(0\le\vartheta\le1\) be smooth, with
\(\vartheta(t_*)=0\) and \(\vartheta(t_0)=1\). For small
\(\varepsilon>0\), set
\[
 v=\Phi+\varepsilon\eta\vartheta.
\]
With \(\vartheta\) fixed, Lemma~\ref{lem:strict-stability}, or
Lemma~\ref{lem:frozen-stability} in the remaining singular range,
gives, for sufficiently small \(\varepsilon\),
\[
 \partial_t\beta(v)+\Lp v\le-\delta
       \quad\hbox{in }D\times(t_*,t_0).
\]
Locally uniform convergence of the symmetric principal values
allows symmetric integration, so the same inequality holds weakly. Moreover, \(v\in C([t_*,t_0];X)\), it is
\(C^{2,1}\) near the comparison cylinder, and it remains positive
there. Since \(\Phi\le u\), one has \(v\le u\) outside
\(\supp\eta\) and on the initial slice. Lemma~\ref{lem:comparison}
gives \(v\le u\) up to time \(t_0\), contrary to
\(v(z_0)=u(z_0)+\varepsilon\).

For a test from above with a positive residual, replace the
interior profile by
\(\phi(y,t_0)+(a-\nu)(t-t_0)\). It lies above the original
test at backward times. Subtracting
\(\varepsilon\eta\vartheta\) gives a strict weak supersolution
that lies above \(u\) outside \(\supp\eta\) and initially.
Comparison contradicts its value \(u(z_0)-\varepsilon\) at
\(z_0\). This proves both viscosity inequalities.
\end{proof}

\begin{proof}[Proof of Theorem~\ref{thm:main}]
Apply Proposition~\ref{prop:viscweak} to the subsolution and
supersolution inequalities, and use Proposition~\ref{prop:weakvisc}
for the converse. The cylinders were arbitrary, so both implications hold
throughout \(Q\).
\end{proof}

\section{Nonnegative solutions}
\label{sec:nonnegative}

We now prove Theorem~\ref{thm:nonnegative}. Positive truncations
provide local energy estimates that are uniform as the truncation
level tends to zero. Global nonnegativity then allows us to pass
across zero times by a zero-slice argument. For the converse, we
combine the weak Harnack inequality with the weak equation.
Throughout this section, when \(1<p<2\), viscosity tests are taken
in the variable \(\beta(u)\), as in
Definition~\ref{def:viscosity-singular}.

\subsection{Truncation and local energy}

\begin{proposition}
\label{prop:nonnegative-subsolution}
Let \(u\) satisfy \eqref{eq:admissible-class} and be a viscosity
subsolution in \(Q\). Then
\(u_+\in L^p_\loc(I;W^{s,p}_\loc(\Omega))\) and \(u_+\) is a
local weak subsolution. In particular, if \(u(\cdot,t)\ge0\)
almost everywhere in \(\R^N\) for every \(t\in I\), then
\(u\) is a local weak subsolution.
\end{proposition}

\begin{proof}
For \(0<\varepsilon<1\), define
\[
 u_\varepsilon=\max\{u,\varepsilon\}
       \quad\hbox{on }\R^N\times I.
\]
Since the truncation map is Lipschitz,
\(u_\varepsilon\in C(Q)\cap C_\loc(I;X)\). On every compact
time interval \(J\Subset I\),
\[
 \sup_{t\in J}\int_{\R^N}|u_\varepsilon(x,t)|^{p-1}\varpi(x)\dd x
 \le \sup_{t\in J}\int_{\R^N}|u(x,t)|^{p-1}\varpi(x)\dd x
       +\int_{\R^N}\varpi(x)\dd x.
\]
We claim that \(u_\varepsilon\) is a viscosity subsolution.
Let \(\phi\) touch it from above at \(z_0\) on a backward
cylinder with spatial radius \(\rho\). For \(1<p<2\), we use
the admissible tests in the \(u\)-variable formulation of
Lemma~\ref{lem:positive-tests}.
If \(u(z_0)>\varepsilon\), then \(\phi\) also touches \(u\)
from above. When \(p\ge2\), apply its viscosity subsolution
inequality. When \(1<p<2\), the function \(\psi=\beta(\phi)\)
is \(C^{2,1}\), since \(\phi\ge u_\varepsilon\ge\varepsilon\)
on the whole contact cylinder. It touches \(\beta(u)\) from
above, and \(\gamma(\psi)=\phi\) has the required spatial
admissibility. Definition~\ref{def:viscosity-singular} therefore
gives the same inequality
\(\beta'(u(z_0))\phi_t(z_0)+\Lp\phi_{\rho,u}(z_0)\le0\).
The two spatial patches agree inside the ball, and
\(u_\varepsilon\ge u\) outside it. Since \(J_p\) is increasing,
\[
 \Lp\phi_{\rho,u_\varepsilon}(z_0)
 \le\Lp\phi_{\rho,u}(z_0),
\]
which proves the claim at such points. If
\(u_\varepsilon(z_0)=\varepsilon\), backward contact gives
\(\phi(x_0,t)\ge\varepsilon=\phi(z_0)\) for earlier times,
hence \(\phi_t(z_0)\le0\). The spatial patch is at least
\(\varepsilon\) almost everywhere at time \(t_0\), so its operator
at \(z_0\) is nonpositive. The subsolution inequality follows.

Theorem~\ref{thm:main} implies that each \(u_\varepsilon\) is a weak
subsolution. To obtain energy bounds uniform in \(\varepsilon\),
fix \(D\times J\Subset Q\),
\(0\le\zeta\le1\), \(\zeta\in C_c^\infty(D)\), and
\(0\le\chi\in C_c^\infty(J)\). Apply
Proposition~\ref{prop:energy} to the upper envelopes of
\(u_\varepsilon\), using a common local bound
\(K=2(1+\|u_+\|_{L^\infty(D\times J)})\) on a slightly
smaller cylinder containing the supports of the cutoffs. The
uniform tail bound above controls the separated interactions.
For each fixed \(\varepsilon\), let the envelope parameter tend
to zero. The residual error vanishes, and Fatou's lemma applied
to the subsolution version of \eqref{eq:caccioppoli}, together
with the cutoff seminorm inequality, gives
\begin{equation}
 \sup_{0<\varepsilon<1}\int_I\chi(t)
       [u_\varepsilon(\cdot,t)\zeta]_{s,p}^p\dd t
 \le C,
 \label{eq:floor-energy}
\end{equation}
where \(C\) depends only on \(N,s,p\), the cylinder and cutoffs,
\(K\), and \(\sup_{t\in J}\|u(\cdot,t)\|_X\). In particular,
it is independent of the lower bound
\(u_\varepsilon\ge\varepsilon\).

Since \(0\le u_\varepsilon-u_+\le\varepsilon\), the convergence
is locally uniform in \(Q\), and
\[
 \sup_{t\in J}\int_{\R^N}|u_\varepsilon(x,t)-u_+(x,t)|^{p-1}
          \varpi(x)\dd x
 \le\varepsilon^{p-1}\int_{\R^N}\varpi(x)\dd x.
\]
Thus \(u_\varepsilon\to u_+\) also in the tail topology,
uniformly on compact time intervals. Estimate~\eqref{eq:floor-energy}
and Fatou's lemma give the local energy membership of \(u_+\).
For the separated exterior interactions,
the bound
\(|u_\varepsilon|^{p-1}\le |u|^{p-1}+1\)
provides an integrable majorant through the weighted tail
assumption.
The flux argument in Proposition~\ref{prop:viscweak} then passes
to the limit in the weak subsolution inequalities.
\end{proof}

\subsection{Zero slices and the weak formulation}

For a globally nonnegative viscosity supersolution, the constant
zero test forces an entire spatial slice to vanish whenever it
vanishes at one interior point.

\begin{lemma}
\label{lem:viscosity-zero-slice}
Let \(u\) satisfy \eqref{eq:admissible-class} and be a viscosity
supersolution in \(Q\), with \(u(\cdot,t)\ge0\) almost
everywhere in \(\R^N\) for each \(t\in I\). If
\(u(x_0,t_0)=0\) at an interior point, then
\(u(\cdot,t_0)=0\) in \(X\), and hence pointwise in \(\Omega\).
\end{lemma}

\begin{proof}
For every sufficiently small \(\rho>0\), use the zero test on
a backward cylinder at \((x_0,t_0)\). If \(p\ge2\), this
means \(\phi\equiv0\) as a test for \(u\). If \(1<p<2\),
use \(\psi\equiv0\) as a test for \(\beta(u)\); this test is
admitted at zero by Definition~\ref{def:viscosity-singular}.
In either case the viscosity inequality gives
\[
 0\le-2\int_{\R^N\setminus B_\rho(x_0)}
       \frac{u(y,t_0)^{p-1}}{|x_0-y|^{N+sp}}\dd y\le0.
\]
Thus \(u(\cdot,t_0)=0\) almost everywhere outside this ball.
Letting \(\rho\downarrow0\) proves equality to zero in \(X\).
Continuity in \(\Omega\) gives the pointwise assertion there.
\end{proof}

\begin{proposition}
\label{prop:nonnegative-viscweak}
A viscosity solution satisfying \eqref{eq:admissible-class},
with \(u(\cdot,t)\ge0\) almost everywhere in \(\R^N\) for
every \(t\in I\), is a local weak solution in \(Q\).
\end{proposition}

\begin{proof}
The local energy membership follows from
Proposition~\ref{prop:nonnegative-subsolution}. Fix \(x_*\in\Omega\) and put
\[
 P=\{t\in I:u(x_*,t)>0\},\qquad Z=I\setminus P.
\]
The set \(P\) is open. Lemma~\ref{lem:viscosity-zero-slice}
shows that \(u>0\) on \(\Omega\times P\) and
\(u(\cdot,t)=0\) in \(X\) for every \(t\in Z\).
By Theorem~\ref{thm:main}, the weak equation holds on every
interval component of \(P\). It remains to join these identities
across \(Z\).

Let \(\eta\in C_c^\infty(Q)\), and choose \([a,b]\Subset I\)
so that \(\eta\) vanishes near both time endpoints and is
supported in \(\Omega\times(a,b)\). Set
\[
 f(t)=\int_\Omega\beta(u(x,t))\eta(x,t)\dd x,
 \qquad
 g(t)=\int_\Omega\beta(u(x,t))\eta_t(x,t)\dd x
       -\E(u(\cdot,t),\eta(\cdot,t)).
\]
We have \(f\in C([a,b])\), while local energy and the uniform
\(X\)-bound give \(g\in L^1(a,b)\) by
Lemma~\ref{lem:energy-testing}. On each component \((c,d)\) of
\(P\cap(a,b)\), the weak equation gives \(f'=g\) in
distributions. Thus \(f\) is locally absolutely continuous there.
Integrability of \(g\) and continuity of the endpoint traces yield
\[
 \int_c^d g(t)\dd t=f(d)-f(c)=0.
\]
Indeed, each endpoint either belongs to \(Z\), where
\(\beta(u)=0\), or is \(a\) or \(b\), where \(\eta=0\).
On \(Z\cap(a,b)\), the full spatial slice vanishes almost
everywhere, so \(g=0\). Since \(g\in L^1(a,b)\), we may sum
over the countably many components of \(P\cap(a,b)\) and obtain
\(\int_a^b g=0\). This is \eqref{eq:weak}.
\end{proof}

\subsection{Viscosity tests at zeros of weak solutions}

For weak solutions, the zero-slice property follows from the weak
Harnack inequality and the weak equation.

\begin{lemma}
\label{lem:weak-zero-slice}
Let \(u\) satisfy \eqref{eq:admissible-class} and be a local
weak solution in \(Q\). Suppose that \(u(\cdot,t)\ge0\)
almost everywhere in \(\R^N\) for every \(t\in I\). If
\((x_0,t_0)\in Q\) and \(u(x_0,t_0)=0\), then
\(u(\cdot,t_0)=0\) in \(X\). There also exist times
\(t_j<t_0\), \(t_j\to t_0\), such that
\(u(\cdot,t_j)=0\) in \(X\). All these slices vanish
pointwise in \(\Omega\).
\end{lemma}

\begin{proof}
We apply \cite[Theorem~1.2]{Prasad24}. Its assumptions hold for
every \(p>1\) in the present setting: local energy and the
uniform weighted tail bound are available, and local boundedness
implies the required local \(L^\infty_tL^2_x\) bound.
The continuous interior traces give the endpoint weak formulation
by temporal cutoffs, as in Remark~\ref{rem:time-traces};
Lemma~\ref{lem:energy-testing} and density justify the energy
tests. No differentiability of \(\beta\) at zero is needed.

The cited theorem gives a spatial essential-infimum estimate at
every time in a subsequent interval. Taking the essential infimum
over that interval, we obtain constants \(q\in(0,1)\),
\(\delta\in(0,1)\), and \(H>0\), depending only on \(N,s,p\),
such that, on sufficiently small interior cylinders,
\[
 \left(\frac{1}{|B_r|}\int_{B_r(x_0)}u(x,\tau)^q\dd x\right)^{1/q}
 \le H\operatorname*{ess\,inf}_{B_{2r}(x_0)\times
       (\tau+\delta(r/2)^{sp},\,\tau+\delta r^{sp})}u.
\]
Choose \(r\) so small that
\(B_{8r}(x_0)\times[t_0-\delta r^{sp},t_0+(8r)^{sp}]\Subset Q\),
which contains the enlarged cylinders required by the estimate. Set
\[
 J_r=(t_0-\delta r^{sp},\ t_0-\delta(r/2)^{sp}).
\]
For each \(\tau\in J_r\), the point \((x_0,t_0)\) belongs to
the interior of the cylinder on the right. Since \(u\ge0\) and
\(u(x_0,t_0)=0\), continuity makes its essential infimum zero.
The estimate and interior continuity therefore give
\(u=0\) on \(B_r(x_0)\times J_r\).

Take a nonzero \(0\le\zeta\in C_c^\infty(B_r(x_0))\) and
\(\chi\in C_c^\infty(J_r)\). In the weak equation tested with
\(\zeta\chi\), the time term and the interaction within
\(B_r(x_0)\) vanish. Consequently
\[
 0=-2\int_{J_r}\chi(t)\int_{B_r(x_0)}\zeta(x)
       \int_{\R^N\setminus B_r(x_0)}
       \frac{u(y,t)^{p-1}}{|x-y|^{N+sp}}\dd y\dd x\dd t.
\]
To pass from almost every time to every time, set
\[
 A(y)=\mathbf1_{\R^N\setminus B_r(x_0)}(y)
       \int_{B_r(x_0)}\frac{\zeta(x)}{|x-y|^{N+sp}}\dd x.
\]
The weight \(A\) is positive outside the ball, and
\(A(y)\le C\varpi(y)\) since \(\supp\zeta\Subset B_r(x_0)\).
By Lemma~\ref{lem:tail-continuity},
\(t\mapsto\int A(y)u(y,t)^{p-1}\dd y\) is continuous.
The weak identity makes it zero almost everywhere on \(J_r\),
and hence everywhere there. Positivity of \(A\) and the local
vanishing already established imply
\(u(\cdot,t)=0\) in \(X\) for every \(t\in J_r\).
Choose \(r_j\downarrow0\) and \(t_j\in J_{r_j}\). Then
\(t_j\to t_0\), and continuity into \(X\) gives
\(u(\cdot,t_0)=0\) in \(X\). Interior continuity gives the
pointwise conclusion on \(\Omega\).
\end{proof}

\begin{proof}[Proof of Theorem~\ref{thm:nonnegative}]
The viscosity-to-weak implication is
Proposition~\ref{prop:nonnegative-viscweak}, and the one-sided
assertion is Proposition~\ref{prop:nonnegative-subsolution}.
Conversely, suppose that \(u\) is a weak solution. At a point
where \(u>0\), Proposition~\ref{prop:weakvisc} applies locally.
Consider a contact point \(z_0=(x_0,t_0)\) with \(u(z_0)=0\).
By Lemma~\ref{lem:weak-zero-slice}, the full spatial slice at
\(t_0\) vanishes almost everywhere, and pointwise in \(\Omega\).

Suppose first that \(1<p<2\). Put \(v=\beta(u)\) and let
\(\psi\) touch \(v\) from below at \(z_0\). Since \(\psi(z_0)=0\),
identity~\eqref{eq:zero-transformed-operator}
reduces the spatial term to a linear principal value. This is
finite because \(sp<2\). The patch is nonpositive at time
\(t_0\), so the spatial term is nonnegative. The earlier zero slices
from
Lemma~\ref{lem:weak-zero-slice} give
\(\psi(x_0,t_j)\le0=\psi(x_0,t_0)\), hence
\(\psi_t(z_0)\ge0\). This proves the supersolution inequality
in Definition~\ref{def:viscosity-singular}.
For a test from above, the spatial term is nonpositive, while
backward contact and \(v\ge0\) give \(\psi_t(z_0)\le0\).
Thus the subsolution inequality also holds.

It remains to consider \(p\ge2\). If \(\phi\) touches from
below, its spatial patch at \(t_0\)
is nonpositive almost everywhere and equals zero at \(x_0\). Hence
\(\Lp\phi_{\rho,u}(z_0)\ge0\). For \(p>2\) the time
coefficient \(\beta'(0)\) is zero, which proves the
supersolution inequality. For \(p=2\), use the earlier zero
times from Lemma~\ref{lem:weak-zero-slice}: for all sufficiently
large \(j\), \(\phi(x_0,t_j)\le0=\phi(x_0,t_0)\), so
\(\phi_t(z_0)\ge0\). This proves the same inequality.

If \(\phi\) touches from above, its spatial patch at \(t_0\)
is nonnegative almost everywhere, whence
\(\Lp\phi_{\rho,u}(z_0)\le0\).
For \(p>2\) this suffices. For \(p=2\), backward contact and
\(u\ge0\) give \(\phi(x_0,t)\ge0\) for earlier contact
times, and thus \(\phi_t(z_0)\le0\). The viscosity subsolution
inequality follows, completing the proof.
\end{proof}

\begin{corollary}
\label{cor:zero-times}
Let \(u\) be a weak or viscosity solution under the assumptions of
Theorem~\ref{thm:nonnegative}. If \(u(x_0,t_0)=0\) at some point of \(Q\),
then \(u(\cdot,t)=0\) in \(X\) for every \(t\in I\) with \(t\le t_0\).
In particular, once \(u\) is positive at one point of \(\Omega\), it is
strictly positive throughout \(\Omega\) at every later time in \(I\).
\end{corollary}

\begin{proof}
By Theorem~\ref{thm:nonnegative} we may use the weak equation.
If \(u(x_0,t)>0\) for some \(t<t_0\), the component of
\(\{\tau\in I:u(x_0,\tau)>0\}\) containing \(t\) has a right endpoint
\(b\le t_0\), with \(u(x_0,b)=0\). Lemma~\ref{lem:weak-zero-slice}
gives zero slices at times approaching \(b\) from the left, contrary
to positivity on that component. Thus \(u(x_0,t)=0\) for every
\(t\in I\) with \(t\le t_0\). Applying the same lemma at each such
time gives the assertion.
\end{proof}

\section{Sign conditions, exterior data, and viscosity tests}
\label{sec:scope}

\subsection{Sign conditions}
Both \(J_p\) and \(\Lp\) are odd. Applying Theorem~\ref{thm:main} to
\(-u\) gives equivalence on negative cylinders, and hence locally on
\(\{u\ne0\}\) for functions satisfying \eqref{eq:admissible-class}.
Theorem~\ref{thm:nonnegative} likewise applies, after a change of sign,
to globally nonpositive functions. When \(p=2\), adding a constant
makes \(u\) positive on any given compact interior cylinder and preserves
both formulations as well as the tail condition, since \(1\in X\).
Thus equivalence holds for arbitrary signed functions in the linear case.

For \(p\ne2\), suppose \(u\ge0\) in \(Q\) and allow negative
exterior values. At an interior zero the integral in
Lemma~\ref{lem:viscosity-zero-slice} may then contain terms of both
signs. Moreover, the truncations \(\max\{u,\varepsilon\}\) used in
Proposition~\ref{prop:nonnegative-subsolution} recover the exterior
values of \(u_+\). The argument thus leaves open equivalence at such
a zero for the prescribed exterior data.

\subsection{Exterior data and tail continuity}
The tail-space assumption permits spatially discontinuous and
unbounded exterior values. For example, suppose that
\(u\in C(Q)\), that \(t\mapsto u(x,t)\) is continuous for
almost every \(x\in\R^N\), and that, for every compact
\(J\Subset I\),
\[
 |u(x,t)|\le C_J(1+|x|)^\gamma
 \quad\hbox{for almost every }x\in\R^N\hbox{ and all }t\in J,
 \qquad 0\le\gamma<\frac{sp}{p-1}.
\]
Then dominated convergence gives \(u\in C_\loc(I;X)\), since
\((1+|x|)^{\gamma(p-1)}\varpi(x)\) is integrable.

Full-space supremal regularization need not preserve \(X\), even
for continuous functions. To see this, let \(q=p-1\), choose
\(0\le\psi\in C_c^\infty(B_1)\) with \(\psi(0)=1\), and put
\(x_k=3^k e_1\), \(\varepsilon_k=2^{-k}\). Define
\[
 A_k=(1+|x_k|^{N+sp})^{1/q},\qquad
 f(x)=\sum_{k=1}^\infty
       A_k\psi\!\left(\frac{x-x_k}{\varepsilon_k}\right).
\]
The supports are disjoint and locally finite, and
\(\|f\|_X^q\le C\sum_k\varepsilon_k^N<\infty\).
Fix \(r>0\). For \(x\in B_{r/2}(x_k)\), taking \(h=x_k-x\) gives
\[
 \sup_{|h|\le r} e^{-\lambda_r|h|^\kappa}f(x+h)
 \ge e^{-\lambda_r(r/2)^\kappa}A_k .
\]
These balls are eventually disjoint, and
\[
 \int_{B_{r/2}(x_k)} A_k^q\varpi(x)\dd x\ge c_r>0
\]
with \(c_r\) independent of \(k\). Summing over \(k\) shows that the
full-space envelope has infinite tail moment. In Lemma~\ref{lem:preserve},
localization confines
the regularization error to the exterior tail. A kernel depending also
on position or time would introduce a further error from translating
the kernel itself.

Changing the exterior values at a single time preserves the weak
identity and can change the viscosity inequality at that time.
Let \(\Omega=B_1\), \(t_0\in I\), and choose
a nonzero \(0\le g\in C_c^\infty(B_3\setminus\overline{B_2})\).
Define
\[
 u(x,t)=1+g(x)\mathbf1_{\{t=t_0\}}.
\]
This function is continuous in \(Q\), is globally bounded and
strictly positive, and agrees almost everywhere in space-time
with the constant weak solution \(1\). It therefore satisfies
the weak equation. For \(p>p_*\), the lower test \(\phi\equiv1\) at
\((0,t_0)\), with any patching radius \(\rho<1\), gives
\[
 \beta'(1)\phi_t(0,t_0)+\Lp\phi_{\rho,u}(0,t_0)
 =-2\int_{\R^N}\frac{g(y)^{p-1}}{|y|^{N+sp}}\dd y<0.
\]
For \(p\le p_*\), use instead
\(\phi_\varepsilon(x,t)=1-\varepsilon|x|^\kappa\), with
\(\kappa>sp/(p-1)\) and sufficiently small \(\varepsilon>0\).
This is an admissible positive lower test. Its interior contribution
is \(C\varepsilon^{p-1}\rho^{\kappa(p-1)-sp}\), whereas the
displayed negative exterior contribution is unchanged. Hence the
supersolution inequality fails for small \(\varepsilon\).
Continuity into \(X\) fixes the tail representative at each time.

\subsection{Viscosity tests at critical points and zeros}
For the spatial profile \(f(x)=1+|x|^2\) near the origin,
the absolute local contribution to its operator is a positive
multiple of
\[
 \int_0^\rho r^{2(p-1)-sp-1}\dd r.
\]
It diverges when \(p\le p_*\), with no cancellation between
opposite directions. This accounts for the flatness condition
\eqref{eq:critical-admissibility} at a critical point of a test.

For \(1<p<2\) and \(u\ne0\), the formal change from
\(\partial_t J_p(u)\) to \(u_t\) amounts to multiplying
\eqref{eq:main} by \(|u|^{2-p}/(p-1)\). Let \(\Omega=B_1\) and let
\(g\ge0\) be a nonzero smooth function supported in
\(B_4\setminus\overline{B_3}\). The time-independent function
\(u(x,t)=g(x)\) vanishes in a neighborhood of \(\Omega\) and
classically satisfies
\[
 u_t+\frac{|u|^{2-p}}{p-1}\Lp u=0\quad\hbox{in }Q.
\]
For the original equation the residual is
\begin{equation}
 \partial_t\beta(u)+\Lp u
 =-2\int_{\R^N}\frac{g(y)^{p-1}}{|x-y|^{N+sp}}\dd y<0
 \quad\hbox{in }Q.
 \label{eq:exterior-bump-residual}
\end{equation}
The constant zero lower test for \(v=\beta(u)\) in
Definition~\ref{def:viscosity-singular} detects this same exterior term.

\section*{Statements and declarations}

\subsection*{Funding}
This work was supported by the National Natural Science Foundation
of China (No.~12471128).

\subsection*{Competing interests}
The authors declare no competing interests.

\subsection*{Data availability}
No datasets were generated or analyzed in this theoretical study.

\subsection*{AI disclosure}
During the preparation of this manuscript, the authors used
ChatGPT-6 Astra as an auxiliary tool for editorial assistance and preliminary
exploration of proof ideas. All computations, proofs, references, and conclusions
were independently verified by the authors, who take full responsibility for the
content of this paper.

\end{document}